\documentclass{article}%
\usepackage{amsfonts}
\usepackage{amsmath,amssymb}
\usepackage[applemac]{inputenc}
\usepackage{amsmath,amssymb,fullpage}
\usepackage{color}
\usepackage{amsmath}
\usepackage{amssymb}
\usepackage{graphicx}
\usepackage[margin=1.5in]{geometry}
\providecommand{\U}[1]{\protect\rule{.1in}{.1in}}

\newcommand{\dproof}{\noindent {Proof.} \quad}
\newcommand{\fproof}{\hfill $\square$ \bigskip}
\newtheorem{definition}{Definition}[section]
\newtheorem{theorem}[definition]{Theorem}

\newtheorem{remark}[definition]{ \it Remark}

\numberwithin{equation}{section}

\def\1B{\text{1\!\!I}}

\newcommand{\mc}{\mathcal}

\def\mb{\mathbb}
\def\mc{\mathcal}

\title{Stochastic optimal control of  pre-exposure prophylaxis for HIV infection}
\author{Kristina Rognlien Dahl\footnote{kristrd@math.uio.no -
Department of Mathematics, University of Oslo, P.O. Box 1053 Blindern, N-0316 Oslo, Norway}
		\and 
		Jasmina Djordevi\'c\footnote{nina19@pmf.ni.ac.rs, djordjevichristina@gmail.com - Faculty of Science and Mathematics, University of Ni\v{s},
Vi\v{s}egradska 33, 18000 Ni\v{s}, Serbia}}

 \definecolor{mypink1}{rgb}{0.858, 0.188, 0.478}
\begin{document}
\maketitle

\begin{abstract}

In this paper, we study the stochastic optimal control problem for the PReP vaccine in the stochastic model for HIV/AIDS with PReP. By using the stochastic maximum principle, we derive the stochastic optimal control of PReP for the unconstrained control problem, as well as for two different types of budget constrains. We illustrate the results by a numerical example. We first study the PReP stochastic differential equation dynamics with a constant, deterministic PReP treatment rate. Then, we compare this to the stochastic optimal control in the unconstrained case.
\bigskip


\medskip


\end{abstract}

\section{Introduction}

The HIV virus continues to be a major global public health issue, still taking millions of lives. Effective antiretroviral drugs can control the virus and help prevent transmission so that people with HIV, and those at substantial risk, can enjoy long and healthy lives (WHO 2015). One successful medicine is antiretroviral therapy (ART), which has shown globally positive results. 

Recently, the World Health Organization's Global Health Sector Strategy on HIV recommended that those at substantial risk of HIV infection should be offered pre-exposure prophylaxis (PReP) as a prevention measure for the reduction of new HIV  infections. PReP is an antiretroviral medication to prevent the acquisition of HIV infection by uninfected persons. It is considered an effective and safe mechanism for preventing  HIV infections (WHO 2015). Only people who are HIV-negative and at very high risk for HIV infection can be under PReP treatment. 

In the literature, there exists several deterministic models based on systems of ordinary differential equations for describing the spread of HIV virus under PReP, see e.g. Silvia and Torres \cite{SilviaTorres} and Campos et. al \cite{Campos}. In this paper, we instead consider a stochastic model for spread of HIV under PReP. Djordjevi\'c and Silva \cite{DjordevicSilva} recently introduced a stochastic model for the spread of HIV/AIDS under PReP treatment according to the following system of stochastic differential equations (SDEs),

\begin{equation}
\label{eq: SDE_system-add}
\begin{array}{llll}
dS(t) &=& [\Lambda {-} \beta (I(t) + \eta_C C(t) + \eta_A A(t))S(t) - \mu S(t) - {\psi}S(t)  + \theta E(t)]dt \\[\smallskipamount]
&&-\sigma(I(t) + \eta_C C(t) + \eta_A A(t))S(t)dB(t) \\[\medskipamount]
dI(t) &=& [\beta (I(t) + \eta_C C(t) + \eta_A A(t))S(t)- \xi_3 I(t) + \alpha A(t) {+} w C(t)] dt \\[\smallskipamount]
&& + \sigma(I(t) + \eta_C C(t) + \eta_A A(t))S(t)dB(t) \\[\medskipamount]
dC(t) &=& [\phi I(t) - \xi_2 C(t)]dt \\[\medskipamount]
dA(t) &=& [\rho I(t) - \xi_1 A(t)]dt \\[\medskipamount]
dE(t) &=& [{\psi }S(t) - \xi_4 E(t)]dt,
\end{array}
\end{equation}
where $B(t)$ is standard Brownian motion, and the variables of the model are  $S$ - susceptible individuals, $I$ - infected, $C$ - chronic stage, $A$ - with AIDS clinical symptoms, $E$ - under PReP. $\Lambda$ is a constant rate under which susceptible individuals increase, it is assumed to have value 2.1$\mu$, where the constant rate $\mu$ is a death rate (with value 1/69.54, Silva and Torres  \cite{SilvaTorres2017}). The model describes how susceptible individuals can get infected through contact with individuals from classes $I, C$ and $A,$ according to the force of infection

\begin{equation}
\label{eq: force_of_infection}
(\beta+\sigma B(t))(I+\eta_C C+\eta_A A),
\end{equation}

\noindent where $\beta$ is the effective contact rate for HIV transmission. Here, $\eta_A$ accounts for the relative infectiousness of individuals with AIDS symptoms belonging to category $A$ (this has value 1.3 according to Silva and Torres \cite{SilvaTorres2017}). Furthermore, $\eta_C$ is the partial restoration of  immune function of individuals with HIV infection that use ART correctly (takes value 0.015, Silva and Torres \cite{SilvaTorres2017}).

Individuals under PReP are transferred to the class $E$ at a rate $\psi$ (this parameter takes value 0.1, Nichols et al.  \cite{Nichols} and Silva and Torres \cite{SilvaTorres}). On the other hand, individuals who stop taking PReP return to the class $S$, at a rate $\theta$ (takes value 0.001 Silva and Torres \cite{SilvaTorres}). Individuals in the infected group, $I$, progress to the class of individuals with HIV infection under ART, $C$, at a rate $\phi$ (takes value 1, Perelson  \cite{Perelson} and Silva and Torres   \cite{SilvaTorres}). The individuals who do not take ART progress to the AIDS class $A$, at rate $\rho$ (takes value 0.1, Sharomi et al. \cite{Sharomi} and Silva and Torres  \cite{SilvaTorres}). The chronic group, $C$, increases with a  rate $\phi$ with the entry of individuals from the class $I$ that are under ART and decreases at a rate $\omega +\mu$ due to the absence  of ART and natural death (it is assumed that $\omega$ takes value 0.09). The evolution of the individuals with  AIDS symptoms is given by the entrance of HIV-infected  individuals that stop ART, at a rate $\rho$, and absence of the individuals that suffer from an AIDS induced death, with the  rate $d$ (takes value 1, Zwahlen and Egger  \cite{Zwahlen}), and natural death, with the constant rate $\mu$. HIV-infected individuals with AIDS symptoms $A$, move to the class of HIV-infected individuals $I$, at a rate $\alpha$.   A fraction $\psi$, where $\psi\in[0,1]$, of susceptible individuals have access to PReP and are transferred to the class $E$. The individuals that stop PReP become susceptible individuals again, at a rate $\theta$, and are transferred to the class $S$. Also, in equation \eqref{eq: SDE_system-add}, $\xi_1 =\alpha+\mu+d,\xi_2 =\omega+\mu, \xi_3=\rho+\phi+\mu$ and  $\xi_4=\mu+\theta$.

Mathematical modeling of processes in biology and medicine, in particular in epidemiology, has led to significant scientific advances both in mathematics and biosciences in areas of prediction and control. The reason for considering stochastic models as opposed to deterministic ones, is that since stochastic models capture randomness they pose a more realistic model of natural events than deterministic ones. The argument against stochastic models is that they are more complicated, and that numerical solutions may be slower than deterministic ones. However, for the PReP model considered in this paper, the stochastic control problem is solvable, and the computational time for the numerical example is small.

Applications of mathematics in biology are opening new pathways of interactions. This is in particular true in the area of (stochastic) optimal control: a branch of applied mathematics that deals with finding control laws for dynamical systems over a period of time such that an objective functional is optimised. In this paper, we consider the controlled stochastic model for spread of HIV with possibility for PReP treatment. The capacity of the available PReP vaccine is limited due to either their costs of production, transport etc. The aim is to use stochastic optimal control theory to determine the optimal percentage of susceptible individuals to be exposed to the PReP vaccine at each time. More about stochastic optimal control theory and the set up of our control problem will be presented in the following sections. 
 
 The rest of the paper is organised as follows; In Section 2, the  stochastic optimal control problem for the PReP problem is introduced, and the existence and uniqueness of the global positive solution for the introduced  control system is  proven. In Section 3, the unconstrained stochastic optimal control of PReP problem is defined, and its solution via Hamiltonian techniques is described.  Section 4 is dedicated to the generalised Lagrange multiplier methods for stochastic optimal control, where two types of constrains for the control problem are introduced and complete proofs for its solutions in most general cases are given. In Section 5, results from Section 4 are applied on optimal control of PReP with budget constraint. In Section 6 some conclusion marks are given and ideas for future research. Section 7 is dedicated to the Appendix to the proofs. At the end of the paper, literature is listed.

\section{The stochastic control model}

In this section, we study a numerical example of the solution of the PReP SDE \eqref{eq: SDE_system-add}, and then introduce the stochastic optimal control PReP problem. We will also show existence of a unique solution to this control SDE.

\subsection{The PReP SDE model: A numerical example} 
\label{sec: prep_sde}

Before we introduce the PReP stochastic optimal control model, we will study equation \eqref{eq: SDE_system-add} for specific choices of $\psi$. This is done in order to get an overview of the dynamics of equation \eqref{eq: SDE_system-add} and the effect of PReP. 

An overview of the choice of values for the parameters in the model is shown in Table \ref{fig: table}.

\begin{figure}
\begin{tabular}{|lll|l|}
  \hline
  Symbol & Description & Value & Reference \\ \hline
  $N = N(0)$ & Initial population & 10 200 &  Assumed based on initial conditions \\ 
  $\mu$ &  Natural death rate & $1/69.54$ &  Silva and Torres \cite{SilvaTorres}\\
  $\Lambda$ & Recruitment rate & $N*\mu$ & Silva and Torres  \cite{SilvaTorres}\\
  $\tilde{\beta}$ & HIV transmission rate & 0.752 & Silva and Torres  \cite{SilvaTorres}\\
  $\beta$ & Scaled HIV transmission rate & 0.752/$N$ & Silva and Torres \cite{SilvaTorres} \\
  $\eta_A$ & Modification parameter & 1.35 & Silva and Torres  \cite{SilvaTorres}\\
  $\eta_C$ & Modification parameter & 0.04 & Silva and Torres  \cite{SilvaTorres}\\
  $ \phi$ & HIV treatment rate for $I$ individuals & 1 &  Silva and Torres \cite{SilvaTorres}\\  
  $\rho$ & Default treatment rate for $I$ individuals & 0.1 &  Silva and Torres \cite{SilvaTorres}\\
  $\alpha$ & AIDS treatment rate & 0.33 & Silva and Torres  \cite{SilvaTorres}\\
  $\omega$ & Default treatment rate for $C$ individuals & 0.09 & Silva and Torres  \cite{SilvaTorres}\\
  $d$  & AIDS induced death rate & 1 & Silva and Torres  \cite{SilvaTorres}\\
  $\psi$  & PreP treatment rate & 0 &  Assumed for optimal control \\
  $\theta$  & PreP default rate & 0.001 & Silva and Torres  \cite{SilvaTorres}  \\
  $\sigma$  &  Force of infection noise parameter & 0.2/$N$ & Assumed  \\
  $w_1$  & Infected weight in performance function & 20 & Assumed  \\
  $w_2$  & PreP/cost weight in performance function & 0.3*N  & Assumed  \\
  \hline
\end{tabular}
\label{fig: table}
\end{figure}

Note that the HIV transmission rate, $\beta$, and the force of infection noise parameter, $\sigma$, are scaled w.r.t. the initial population. This is done since the transmission of HIV is dependent on the total number of people in the population. In this initial example, we will assume that 

\[
N=N(0) = S(0) + I(0) + A(0) + C(0) + E(0).
\]

That is, the initial population equals the sum of the initial number of susceptible, infected, AIDS, chronic and those under PreP. We choose $S(0)=10 000$, I(0)=200 and $A(0)=C(0)=E(0)=0$. We let the terminal time $T=25$ years. We solve the SDE \eqref{eq: SDE_system-add} numerically via the stochastic Euler method which is known to be strongly convergent with order $0.5$. For the following simulation, we used time step size $\Delta t = 1/1000$. 

\begin{figure}
\begin{centering}
\includegraphics[width=0.75\textwidth]{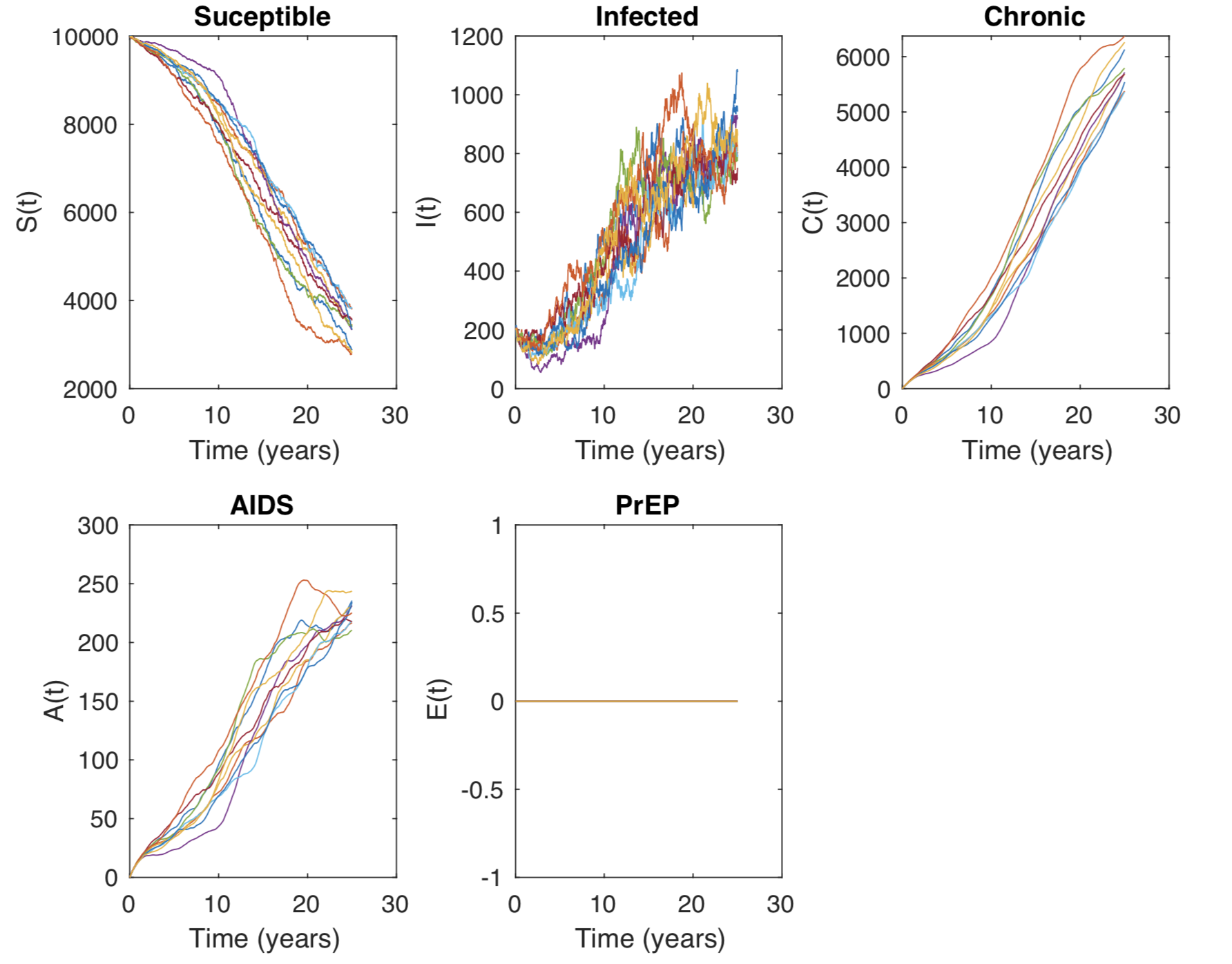}
\end{centering}
\caption{Plot of 10 paths of the solution of the SDE PReP model \eqref{eq: SDE_system-add} over 25 years with no PReP.}
\label{fig: SDE_noPrep}
\end{figure}

In Figure \ref{fig: SDE_noPrep}, we have plotted 10 paths of the solution of the PReP SDE model \eqref{eq: SDE_system-add} with $\psi = 0$. That is, no individuals get PReP treatment. The remaining parameters of the models are chosen as in Table \ref{fig: table}. Note that there is no optimal control involved at this point. Figure \ref{fig: SDE_noPrep} shows the effect of a constant, deterministic PReP treatment rate of 10\% of the group of susceptible individuals. As we can see in Figure \ref{fig: SDE_noPrep}, the number of susceptible individuals decreases, but the number of infected, chronic and AIDS individuals rapidly increases. Note also that the noise in the model is most prominent in the infected and AIDS category. 

\begin{figure}
\begin{centering}
\includegraphics[width=0.75\textwidth]{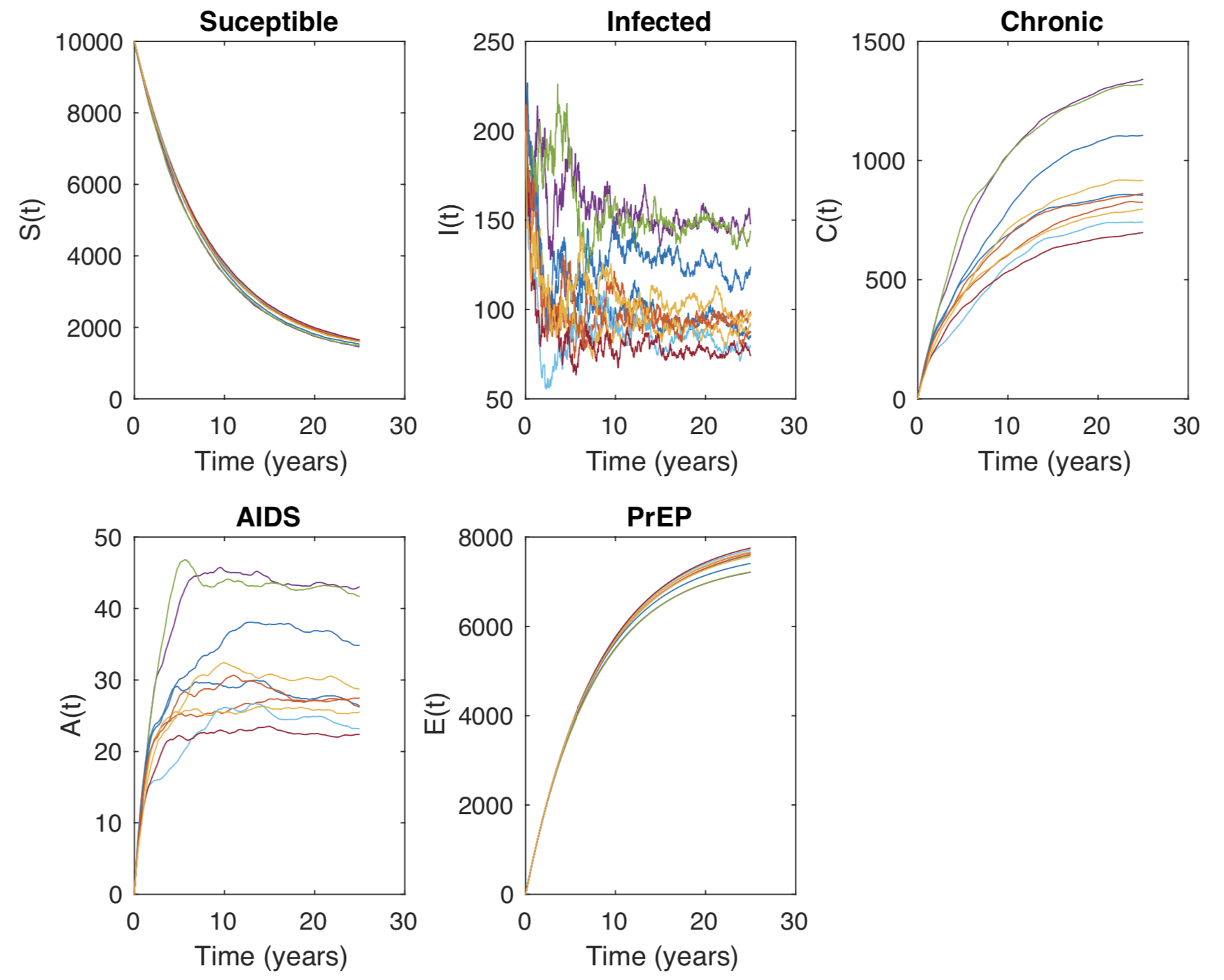}
\end{centering}
\caption{Plot of 10 paths of the solution of the SDE PReP model \eqref{eq: SDE_system-add} over 25 years with $\psi=0.1$, i.e., 10\% of the susceptible individuals get PReP treatment.}
\label{fig: SDE_prep_0point1}
\end{figure}

In Figure \ref{fig: SDE_prep_0point1}, we have plotted 10 paths of the solution of the PReP SDE model \eqref{eq: SDE_system-add} with $\psi = 0.1$. That is, 10\% of all susceptible individuals get PReP treatment. The remaining parameters of the models are chosen as in Table \ref{fig: table}. In Figure \ref{fig: SDE_prep_0point1}, we note that the number of susceptible individuals decreases, and so does the number of infected. However, the decrease in the infected-group is slow and with a large variance. Furthermore, the chronic and AIDS groups increase before stabilizing, but again, the variance is large compared to the no PReP case in Figure \ref{fig: SDE_noPrep}.

\begin{figure}
\begin{centering}
\includegraphics[width=0.75\textwidth]{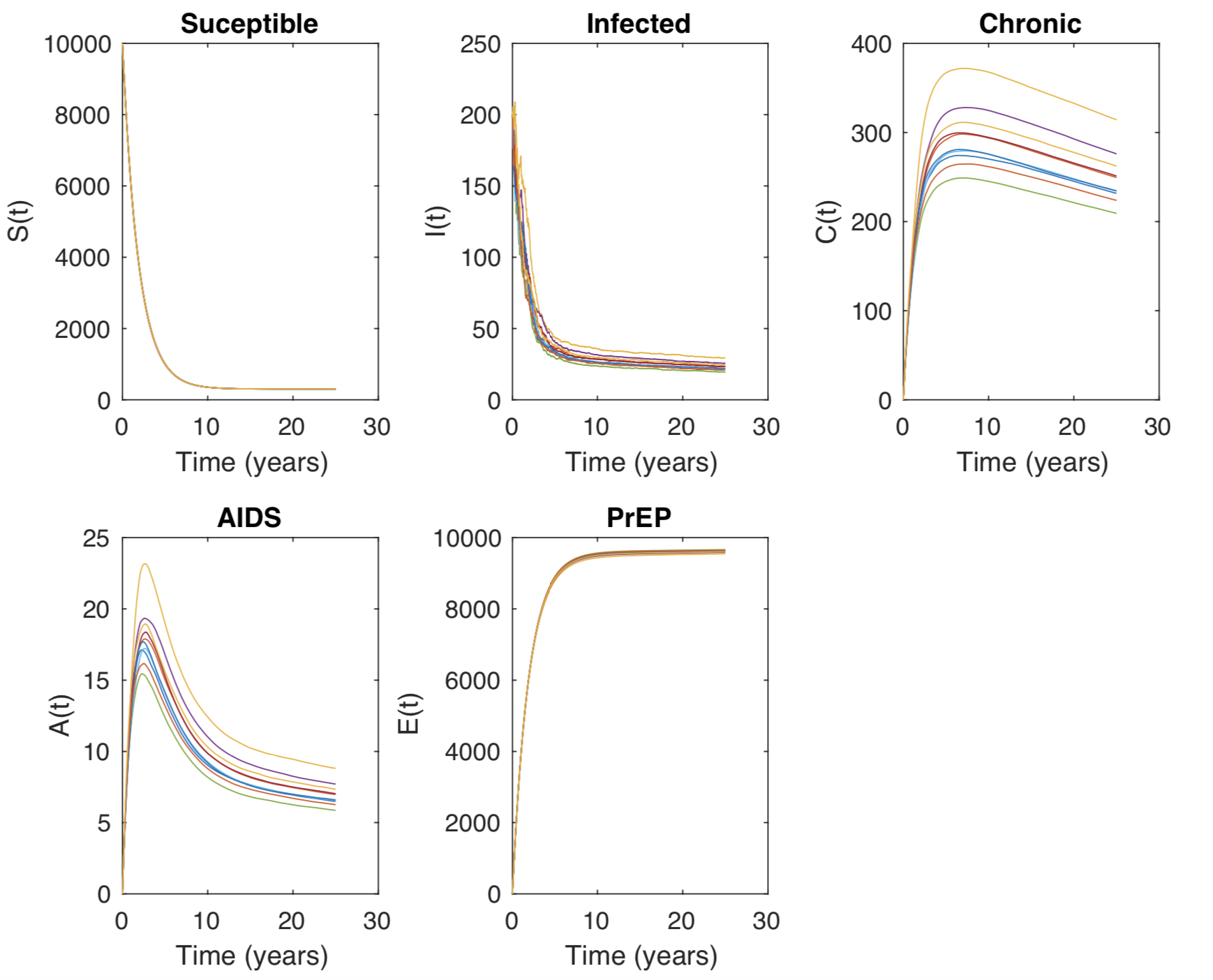}
\end{centering}
\caption{Plot of 10 paths of the solution of the SDE PReP model \eqref{eq: SDE_system-add} over 25 years with $\psi =0.5$, i.e.,  50\% of susceptible individuals get PReP treatment.}
\label{fig: SDE_halfprep}
\end{figure}

In Figure \ref{fig: SDE_halfprep} we have plotted 10 paths of the solution of the PReP SDE model \eqref{eq: SDE_system-add} with $\psi = 0.5$. That is, 50\% of all susceptible individuals get PReP treatment. The remaining parameters of the models are chosen as in Table \ref{fig: table}. Here, we see a rapid decrease, before a stabilisation, in the number of susceptible and infected individuals. For the chronic and AIDS groups, there is an initial increase followed by a gradual decrease. The variance in the infected group is very small in comparison to the no PReP case in Figure \ref{fig: SDE_noPrep} and the low-PReP case in Figure \ref{fig: SDE_prep_0point1}. However, there is still some variance in the chronic and AIDS groups.

The development of the susceptible and chronic groups appear less affected by the noise. Also note that the variance of the processes appear to be larger in the middle case where $\psi=0.1$ than in the "extreme" cases $\psi=0$ and $\psi=0.5$. The expectation and variance of the various processes, $S, I, C, A$ and $E$ can be simulated via Monte Carlo techniques. A more detailed analysis of the numerical aspects of this problem is a work in progress, and will be the topic future works.

\subsection{The stochastic control SDE}

In the previous Section \ref{sec: prep_sde}, we considered the solution of the PReP SDE \eqref{eq: SDE_system-add} under constant, deterministic PReP treatment rates $\psi = 0, 0.1, 0.5$. The purpose of the remaining part of the paper is to introduce and solve a stochastic optimal control problem with the PReP SDE as the state dynamics. The aim is to minimise a performance function which is a weighted sum of a term depending on the number of infected individuals and another term depending on the cost of PReP treatment. A similar analysis was done in Silva and Torres \cite{SilvaTorres} for a deterministic dynamic system. The model in Silva and Torres \cite{SilvaTorres} can be generalised to a stochastic dynamic system by including a Brownian motion in a suitable way. We first introduce the stochastic control PReP SDE.

Let $(\Omega, \mathcal{F}, P)$ be a probability space, where $\Omega$ is the scenario space, $\mathcal{F}$ a $\sigma$-algebra and $P$ the probability measure. We consider continuous time  $t \in [0,T]$. Let $\{B(t)\}_{t \in [0,T]}$ be a Brownian motion, and let  $\{\mathcal{F}_t\}_{t \in [0,T]}$ be the filtration generated by this Brownian motion. In the following, by an adapted processes, we mean adapted with respect to this filtration.

Let $T$ be a terminal time (the final time of interest). Further, let $\{u(t)\}_{t \in [0,T]}$ denote  the stochastic control process, where $u(t, \omega)$ denotes the percentage of individuals under PReP treatment at time $t$ in scenario $\omega \in \Omega$, so $u(t) \in [0,1]$ for all $t \in [0,T]$. Let 

\begin{equation}
\mathcal{A} := \{u(\cdot) \in \mathcal{L}^{\infty}([0,T], \Omega) : 0 \leq u(t) \leq 1 \mbox{ } P\mbox{-a.e.}, u(t) \mbox{ adapted}\}.
\end{equation}

Generalising the set up of Djordevi{\'c} and Silva  \cite{DjordevicSilva}, as well as  Silva and Torres \cite{SilvaTorres}, the system of controlled stochastic differential equations  to model the spread of HIV/AIDS is as follows:

\begin{equation}
\label{eq: SDE_system}
\begin{array}{llll}
dS(t) &=& [\Lambda {-} \beta (I(t) + \eta_C C(t) + \eta_A A(t))S(t) - \mu S(t) - {\psi}S(t) - u(t)S(t) + \theta E(t)]dt \\[\smallskipamount]
&&-\sigma(I(t) + \eta_C C(t) + \eta_A A(t))S(t)dB(t) \\[\medskipamount]
dI(t) &=& [\beta (I(t) + \eta_C C(t) + \eta_A A(t))S(t)- \xi_3 I(t) + \alpha A(t) {+} w C(t)] dt \\[\smallskipamount]
&& + \sigma(I(t) + \eta_C C(t) + \eta_A A(t))S(t)dB(t) \\[\medskipamount]
dC(t) &=& [\phi I(t) - \xi_2 C(t)]dt \\[\medskipamount]
dA(t) &=& [\rho I(t) - \xi_1 A(t)]dt \\[\medskipamount]
dE(t) &=& [{\psi }S(t) - \xi_4 E(t) + S(t)u(t)]dt,
\end{array}
\end{equation}
where $\xi_1 =\alpha+\mu+d,\xi_2 =\omega+\mu, \xi_3=\rho+\phi+\mu$ and  $\xi_4=\mu+\theta$.





The proof of existence of global positive solution of system (\ref{eq: SDE_system}) is similar to the one from  Djordevi{\'c} and Silva \cite{DjordevicSilva}, but with added control. In the sequel we will adjust the proof. 

We will use usual notation $R_+^5=\left\{(x_1,x_2,x_3,x_4,x_5) \mbox{ } | \mbox{ } x_i>0, \mbox{ } i=1, \ldots, 5\right\}$.

\begin{theorem}\label{existance} For any initial value $(S(0), I (0), C(0), A(0), E(0)) \in R^5_+$, there is a unique positive solution $(S(t), I (t), C(t), A(t), E(t))$  of system \eqref{eq: SDE_system}  for every $t \geq 0$ and the solution will remain positive with probability one. That is, $(S(t), I(t), C(t), A(t), E(t))\in R^5_+$ for all $t \geq 0$ almost surely. Moreover, for $N(t)= S(t)+I(t)+ C(t)+ A(t)+E(t)$, it follows that
	
	\begin{equation}
	\label{eq:NgoesLam:div:mu}
	N(t) \to \frac{\Lambda}{\mu} \mbox{ as } t \to \infty.
	\end{equation}

\end{theorem}

\dproof  Since $E(t)=N(t)-S(t)-I(t)-C(t)-A(t),$ for every $t\geq 0$, we reduce system \eqref{eq: SDE_system} to following system of four equations.
	\begin{equation}
	\begin{cases}
	d S(t) = \left[ \Lambda - \beta \left( I(t) + \eta_C \, C(t)  + \eta_A  A(t) \right) S(t) - (\mu+u(t)) S(t) - \psi S(t)\right. \\
	\quad \quad \quad\left.+ \theta (N(t)-S(t)-I(t)-C(t)-A(t))\right] dt\\
	\quad \quad \quad - \sigma \left( I(t) + \eta_C \, C(t)
	+ \eta_A  A(t) \right) S(t) dB(t),\\[0.2 cm]
	d I(t) = \left[  \beta \left( I(t) + \eta_C \, C(t)  
	+ \eta_A  A(t) \right) S(t) - \xi_3 I(t) 
	+ \alpha A(t)  + \omega C(t) \right] dt\\
	\quad \quad \quad 
	+ \sigma \left( I(t) + \eta_C \, C(t)  + \eta_A  A(t) \right) S(t) dB(t), \\[0.2 cm]
	d C(t) = \left[  \phi I(t) - \xi_2 C(t) \right] dt,\\[0.2 cm]
	d A(t) = \left[   \rho \, I(t) - \xi_1 A(t) \right] dt.
	\end{cases}	\label{eq:2}
	\end{equation}
	
	If we prove that there exists a unique positive solution $(S(t), I(t),C(t),A(t))$ of system \eqref{eq:2} for $t\geq 0$, it is equivalent to proving existence of a unique positive  solution of system \eqref{eq: SDE_system}. 
	
	 Hence, for given initial conditions 
	$\left( S(0), I(0), C(0), A(0),E(0) \right) \in R^5_+$, we will prove 
	that there exists a unique positive solution of system \eqref{eq: SDE_system}
	for every $t\geq 0$. Because the coefficients of system \eqref{eq:2}
	are locally Lipschitz continuous, there is a unique local solution 
	on $[0,\tau_0)$ for any initial value $(S(0), I(0), C(0),A(0))$, where $\tau_0$ is known
	in the literature as the \emph{explosion time}. It is necessary to prove that the solution 
	is global, i.e., that $\tau_0=+\infty$ almost surely.
	
	Let $k_0\geq 0$ be sufficiently large such that $S(0), I(0), C(0),A(0)$ lie within the interval $[1/k_0, k_0]$. For each integer $k > k_0$, let us define the stopping time	
	$$
	\tau_k = \inf\left\{t\in [0,\tau_0): \min\{S(t),I(t),C(t),A(t)\}\leq \frac{1}{k} \mbox{ or }  
	\max\{S(t),I(t),C(t),A(t)\}\geq k\right\},
	$$
	where  $\inf \emptyset = \infty$. According to the definition, $\tau_k$ is increasing
	as $k \longmapsto +\infty$. Set $\tau_{\infty} = \lim_{k \longmapsto +\infty}\tau_k$, from what follows $\tau_{\infty}\leq \tau_0$ a.s. In order to complete the proof, we need to prove that $\tau_{\infty}=\infty$.
	
Since the infimum of an empty set is $\infty$ and  $\tau_k \leq  \tau_0$, if we prove that $\tau_+ = \infty$ \textcolor{red}{what is $\tau_+$?} a.s., then the proof of our theorem is complete. Indeed, if $\tau_+ = \infty$ a.s., then $\tau_0 = \infty$, which means that $(S(t),I(t),C(t),A(t)) \in R_4^+$ 
	for $t \geq 0$ a.s. .

	Suppose that there exist a pair of constants 
	$T \geq 0$ and $\epsilon \in (0, 1)$ such that
	$$P(\tau_{\infty}\leq T)\geq \epsilon.$$
	
	Then, there exists $k_1\geq k_0$ such that
	\begin{equation}	
	P(\tau_{k}\leq T)\geq \epsilon \mbox{ for all } k\geq k_1.\label{stop}
	\end{equation}

	Besides, for $t \leq \tau_k$, we have that
	$$ N(t)=\left\{ \begin{array}{lll}
	\frac{\Lambda}{\mu}, & \mbox{ if } S(0)+I(0)+ C(0)+ A(0)+E(0)\le \frac{\Lambda}{\mu} ,\\
	S(0)+I(0)+ C(0)+ A(0)+E(0), &\mbox{ if }S(0)+I(0)+ C(0)+ A(0)+E(0)>\frac{\Lambda}{\mu}. \end{array} \right. :=\tilde{N}. $$
	
	Now, define the twice differentiable function $V:R^4_{+}\longmapsto R_+ \cup\{0\}$ in the following way
	$$V(S, I, C,A) = (S - 1 -\log S) + (I - 1 - \log I) + (C - 1 - \log C)+(A-1-\log A).$$
	The function $V$ is nonnegative ($\log x\leq x-1$ for every $x\geq 0$).
	By applying the It\^o formula to the function $V$, we have
	\begin{equation*}
	\begin{split}
	&dV(S, I, C,A)=\left(1-\frac{1}{S}\right)dS(t)+\frac{1}{2S^2(t)}[dS(t)]^2+\left(1-\frac{1}{I}\right)dI(t)+\frac{1}{2I^2(t)}[dI(t)]^2\\
	&\phantom{dV(S, I, C,A)=}+\left(1-\frac{1}{C}\right)dI(C)+\frac{1}{2C^2(t)}[dC(t)]^2+\left(1-\frac{1}{A}\right)dA(t)+\frac{1}{2A^2(t)}[dA(t)]^2.
	\end{split}
	\end{equation*}

	Hence,
	
	$$dV(S, I, C,A)=K(S,I,C,A)dt+\sigma\left( I(t) + \eta_C \, C(t)
	+ \eta_A  A(t) \right)\frac{I(t)-S(t)}{I(t)}dB(t),$$
	where $K : R^4_+ \longmapsto R_+$ is defined by
	
	\begin{equation*}
	\begin{split}
	& K(S,I,C,A)=\left(1-\frac{1}{S}\right)[\Lambda - \beta \left( I(t) + \eta_C \, C(t) + \eta_A  A(t) \right) S(t) - (\mu+u(t)) S(t) - \psi S(t)  \\
	&\phantom{dV(S, I, C,A)=} + \theta (N(t)-S(t)-I(t)-C(t)-A(t))]+\frac{\sigma ^2}{2}\left(I(t) + \eta_C \, C(t)+ \eta_A  A(t)\right)^2\\
	&\phantom{dV(S, I, C,A)=}+\left(1-\frac{1}{I}\right)[\beta \left( I(t) + \eta_C \, C(t)  + \eta_A  A(t) \right) S(t) - \xi_3 I(t) 
	+ \alpha A(t)  + \omega C(t)]\\
	&\phantom{dV(S, I, C,A)=}+\frac{\sigma^2 S(t) ^2}{2I^2(t)} \left( I(t) + \eta_C \, C(t)  + \eta_A  A(t) \right)^2+\left(1-\frac{1}{C}\right)\left[  \phi I(t) - \xi_2 C(t) \right] \\
	&\phantom{dV(S, I, C,A)=}+\left(1-\frac{1}{A}\right) \left[  \rho \, I(t) - \xi_1 A(t) \right]\\
	&\phantom{dV(S, I, C,A)}\leq\Lambda+\theta (N(t)-S(t)-I(t)-C(t)-A(t))+\beta \left( I(t) + \eta_C \, C(t) + \eta_A  A(t) \right)(S(t)+1)\\
	&\phantom{dV(S, I, C,A)=}+\frac{\sigma ^2}{2}(S^2(t)+I^2(t))\left(I(t) + \eta_C \, C(t)+ \eta_A  A(t)\right)^2\\
	&\phantom{dV(S, I, C,A)=}+(\phi +\rho)I(t)+ \alpha A(t) + \omega C(t)+\mu +1+\psi+\xi_1+\xi_2+\xi_3\\
	&\phantom{dV(S, I, C,A)}=\Lambda+m_1+m_2 \tilde{N}+ m_3\sigma^2 \tilde{N}^2:=\overline{N},
	\end{split}
	\end{equation*}
	where $m_1,m_2,m_3$ are generic constants. We have that the expectation is
\begin{equation}
\begin{split}
E\Big(V(S(\tau_k\wedge T), I(\tau_k\wedge T), C(\tau_k\wedge T),A(\tau_k\wedge T)\Big) \nonumber &\leq E\Big(V(S(0), I(0), C(0),A(0)\Big)+\overline{N}T.
\end{split}
\label{estimate1}
\end{equation}

Let $A_k=\{\tau_k\leq T\}$ for $k\geq k_1$, and from \eqref{stop} it follows that $P(A_k)\geq \epsilon.$ Furthermore, for every $\omega\in a_k$, at least one of the variables $S,I,C$ or $A$ is less than or equal $\frac{1}{k}$, or it is greater or equal with $k$. Then, the
	function $V(S(\tau_k), I(\tau_k), C(\tau_k),A(\tau_k))$ is not less than
	$$k-1-\log k \mbox{ or } \frac{1}{k}-1-\log\frac{1}{k},$$
	i.e.,
	$$V(S(\tau_k), I(\tau_k), C(\tau_k),A(\tau_k))\geq \min\left\{k-1-\log k,\frac{1}{k}-1+\log{k}\right\}.$$
	
	From \eqref{stop} and \eqref{estimate1} it follows that
	
	\begin{equation*}
	\begin{split}
	&E\Big(V(S(0), I(0), C(0),A(0)\Big)+\overline{N}T\}\geq\epsilon \min\left\{k-1-\log k,\frac{1}{k}-1+\log{k}\right\},
	\end{split}
	\end{equation*}
	where $I_{A_k}$ denotes the indicator function of the set ${A_k}$. If we let $k\longmapsto +\infty$, we obtain
	$$+\infty>E\Big(V(S(0), I(0), C(0),A(0)\Big)+\overline{N}T=+\infty,$$
	which is a contradiction. Hence, our assumption $P(\tau_{\infty}\leq T)\geq \epsilon$ is wrong, i.e., it follows that $\tau_{\infty}=\infty$ a.s.
	
	\medskip
	
	If we sum all equations from system \eqref{eq: SDE_system}, then 
	\begin{equation*}
	\begin{split}
	&d(S(t)+I(t)+C(t)+A(t)+E(t))
	=\left[\Lambda-\mu S(t) +( \phi-  \rho - \phi - \mu+\rho) I(t)\right.\\
	&\phantom{d(S(t)+I(t)+C(t)+A(t)+E(t))=}\left.
	+ \left(\alpha  - \alpha - \mu - d \right)A(t)
	+ \left(\omega- \omega - \mu\right) C(t)+(\theta-\theta-\mu)E(t)\right]dt\\
	&\Leftrightarrow \frac{d(S(t)+I(t)+C(t)+A(t)+E(t))}{dt}
	=\Lambda-\mu(S(t)+I(t)+C(t)+A(t)+E(t))-d\cdot A(t).
	\end{split}
	\end{equation*}
	Solving the last equation, we obtain that

	$$
	S(t)+I(t)+ C(t)+ A(t)+E(t)=e^{-\mu t}\left[S(0)+I(0)+ C(0)+ A(0)+E(0)
	+\int_0^t(\Lambda -d\cdot A(t))e^{\mu s}ds\right].
	$$
	Applying L'Hospital's rule, it follows that
	\begin{equation*}
	\begin{split}
	&\lim_{t\rightarrow \infty}[S(t)+I(t)+C(t)+A(t)+E(t)]\\
	&\leq\lim_{t\rightarrow \infty}\frac{S(0)+I(0)+C(0)+A(0)+E(0)+\int_0^t \Lambda e^{\mu s}ds}{e^{\mu t}}\\
	&=\frac{\Lambda}{\mu}.
	\end{split}
	\end{equation*}
	 This completes the proof.
	 
	 \begin{remark}
\label{remark:1}
	It should be noted that the set 
	\begin{equation}
	\Gamma^*=\left\{(S(t),I(t),C(t),A(t),E(t)), S(t)>0,I(t)>0,C(t)>0,A(t)>0,E(t)>0,N(t)\leq \frac{\Lambda}{\mu}\right\}.\label{rang}
	\end{equation}
	is a positively invariant set of system \eqref{eq: SDE_system}   on $\Gamma^*$.
\end{remark}

\fproof


Note that the proof of existence and uniqueness of solution of the control SDE in Theorem \ref{existance} holds even without the reduction from 5 to 4 equations. This reduction is done by using the assumption that

\[
N(t) = S(t) + I(t) + C(t) + A(t) + E(t).
\]
This may not always be the case: For the specific case of HIV, it may be reasonable to assume that individuals who are not sexually active, or who are in a monogamous relationship where neither part has HIV, are not susceptible. However, even if 
\[
N(t) >S(t) + I(t) + C(t) + A(t) + E(t),
\] 
\noindent we still get existence and uniqueness of the control SDE by proceeding as in the proof of Theorem \ref{existance} without the reduction of $E(t)$. This involves a little more notation, but the mathematics are the same.

\begin{remark}
\label{remark:add1}
Note that even though the system (\ref{eq: SDE_system}) depends on the control $u(t)$, it is not a part of the solution of the system. For each time $t$, the control $u(t)$ is a percentage of individuals under PReP treatment. The control has bounded values, between 0 to 1, so the system  (\ref{eq: SDE_system}) can be easily bounded with $u(t)$ as well. Because of this, the conditions for extinction and persistence of the disease for the model (\ref{eq: SDE_system})  do not differ significantly compering to the ones proven by Djordjevi\'c and Silva in \cite{DjordevicSilva} for the stochastic model without $u(t)$.
\end{remark}

Now, we are ready to introduce the PReP stochastic optimal control problem.
	
\section{The stochastic optimal control problem}
\label{sec: control_problem}

In the sequel we will assume that $(S(0), I(0), C(0),A(0),E(0)) \in \Gamma^*$.
Our problem is to determine the PReP strategy $\{u(t)\}_{t \in [0,T]} \in \mathcal{A}$ which minimises the performance functional

\begin{equation}
\label{eq: performance}
J(u) := E\left[\int_0^T [w_1 I(t) + w_2 u^2(t)]dt\right].
\end{equation} 
Here, $w_1, w_2$ are weights given to the number of HIV infected and $u(t)$, which represents the percentage of susceptible individuals under PReP. 

\begin{remark}
\label{remark:2}
The choice of $u^2(\cdot)$ in the performance function is made based on custom in the literature. 
\end{remark}

Silva and Torres \cite{SilvaTorres} include a budget constraint in their deterministic version of the optimal control problem (referred to as a mixed state constraint in Silva and Torres \cite{SilvaTorres}),

\begin{equation}
S(t)u(t) \leq \mathcal{V}, \mbox{ } \mathcal{V} \geq 0 \mbox{ for almost all } t \in [0,T].
\end{equation}



\bigskip
This constraint describes that the number of individuals under PReP should be bounded by a constant $\mathcal{V}$ for almost all times. When generalising this constraint to the stochastic case, we have to alter it slightly in order to be able to solve the constrained problem via Lagrange duality techniques. Let $c(t)$, $t \in [0,T]$ be a given cost function for treating one individual with PReP at time $t$. Two types of constraints are considered:

\begin{equation}
\label{eq: type_i}
E\left[\int_0^T S(t)u(t)c(t)dt\right] \leq \mathcal{V}.
\end{equation}

Let this be a Type $I$ constraint. Alternatively, constraints of the type

\begin{equation}
\label{eq: type_ii}
\int_0^T S(t)u(t)c(t)dt \leq \mathcal{V}, \mbox{ } P\mbox{-a.e}.,
\end{equation}
will be called a Type $II$ constraint. Note that if Type $II$ holds, Type $I$ also holds. In this sense, Type $II$ is a stricter constraint than Type $I$.

First, the stochastic optimal control problem without constraints will be solved, then with respect to both types above. In order to do so, we need to prove some results on solutions of stochastic optimal control problems with constraints, and how such problems can be solved by introducing generalised Lagrange multiplies. These results generalise Theorem 2.1 and 2.2 in Dahl and Stokkereit \cite{DahlStokkereit}.

\section{The unconstrained stochastic optimal control of PReP problem}
\label{sec: unconstrained}

To simplify notation, let us introduce vector

\[
\bold{X}(t) = (X_1(t), \ldots, X_5(t)) := (S(t), I(t), C(t), A(t), E(t)).
\]
Also, let the initial state of the system (the initial condition) be $\bold{x}(0)=(x_{1,0}, x_{2,0}, x_{3,0}, x_{4,0}, x_{5,0})$.

The aim is to solve the following problem:

\begin{equation}
\label{eq: prep_problem}
\min_{u \in \mathcal{A}} J_{x_0}(u) = \min_{u \in \mathcal{A}} E^{x_0} \left[\int_0^T w_1 I(t) + w_2 u^2(t)dt\right]
\end{equation}

We rewrite the system of stochastic differential equations \eqref{eq: SDE_system} in matrix form. To do so, let

\[
\bold{K} := (\Lambda, 0, 0, 0, 0)' ,
\]

\bigskip



and
\[
A := \begin{bmatrix}  - \mu -{\psi} & 0 & 0 & 0 & 0 \\
0 & -\xi_3 & w & \alpha & 0 \\
0 & \phi & -\xi_2 & 0 & 0 \\
0 & \rho & 0 & -\xi_1 & 0 \\
{\psi}  & 0 & 0 & 0 & -\xi_4 \end{bmatrix}, \mbox{ } B:= \begin{bmatrix}  -1 & 0 & 0 & 0 & 0 \\
0 & 0 & 0 & 0 & 0 \\
0 & 0 & 0 & 0 & 0 \\
0 & 0 & 0 & 0 & 0 \\
1 & 0 & 0 & 0 & 0 \end{bmatrix}
\]


\[
\begin{array}{lll}
f(\bold{X}(t)) &=& (-\beta (I(t) + \eta_C C(t) + \eta_A A(t))S(t), \beta (I(t) + \eta_C C(t) + \eta_A A(t))S(t), 0, 0, 0)' \\[\smallskipamount]
&:=& (f_1(\bold{X}(t)), f_2(\bold{X}(t)), 0, 0, 0)',
\end{array}
\]

\[
\begin{array}{llll}
g(\bold{X}(t)) &=& (-\sigma (I(t) + \eta_C C(t) + \eta_A A(t))S(t), \sigma (I(t) + \eta_C C(t) + \eta_A A(t))S(t), 0, 0, 0)' \\[\smallskipamount]
&:=& (g_1(\bold{X}(t)), g_2(\bold{X}(t)), 0, 0, 0)'.
\end{array}
\]

Then, we can rewrite the system in matrix form:

\begin{equation}
\label{eq: matrix_control}
\begin{array}{llll}
d\bold{X}(t) &:=& \{\bold{K} + f(\bold{X}(t)) + A\bold{X}(t) + B\bold{X}(t)u(t)\}dt + g(\bold{X}(t))dB(t) \\[\medskipamount]
&:=& b(\bold{X}(t), u(t))dt + \sigma(\bold{X}(t))dB(t),
\end{array}
\end{equation}

where,

\[
\begin{array}{lll}
b(\bold{X}(t), u(t)) &:=& \bold{K} + f(\bold{X}(t)) + A\bold{X}(t) + B\bold{X}(t)u(t) \\[\smallskipamount]
\sigma(\bold{X}(t)) &:=& g(\bold{X}(t)).
\end{array}
\]

From the matrix form \eqref{eq: matrix_control},  the stochastic optimal control problem can be rewritten as:

\begin{equation}
\label{eq: control_final}
\begin{array}{lll}
\min_{u \in \mathcal{A}} J_{x_0}(u) \\[\medskipamount]
\mbox{s.t.} \\[\medskipamount]
d\bold{X}(t)= b(\bold{X}(t), u(t))dt + \sigma(\bold{X}(t))dB(t), \mbox{ } t \in [0,T].
\end{array}
\end{equation}

This is a stochastic optimal control problem of the standard form given by  {\O}ksendal \cite{Oksendal}. From Theorem \ref{existance}  we have that for each $u \in \mathcal{A}$, there exists a unique solution to the controlled SDE \eqref{eq: matrix_control}.




Note that the SDE \eqref{eq: matrix_control} is Markovian. Hence, we can solve the stochastic control problem \eqref{eq: control_final} either via stochastic maximum principles or by using stochastic dynamic programming. Following the deterministic case in Silva and Torres \cite{SilvaTorres}, we choose the maximum principle approach.

To the best of our knowledge, the earliest works on a stochastic maximum principle are Kushner \cite{Kushner}, Bismut \cite{Bismut} and Peng \cite{Peng} who all considered the no-jump case. A necessary maximum principle was derived in the jump case by Tang and Li \cite{TangLi}. A sufficient maximum principle in the jump case is given in Framstad et al. \cite{FOS}. Following this, many variants of the stochastic maximum principle have been derived. Examples include Baghery and {\O}ksendal \cite{Baghery} for partial information, {\O}ksendal and Sulem \cite{OksendalSulem} for delay and Buckdahn et al. \cite{Buckdahn} for mean-field systems. 

\medskip

Our PReP control problem is of the no-jump type, ie it is of Brownian form studied in the earliest works by Kushner \cite{Kushner}, Bismut \cite{Bismut} and Peng \cite{Peng}. See the appendix for a brief summary of the theory of stochastic maximum principle for jump diffusions.

\begin{remark}
\label{remark:3}
It would also be possible to solve this problem via stochastic dynamic programming. In this case, one would derive the Hamilton Jacobi Bellman (HJB) partial differential equation. By solving this PDE, one can derive the optimal value function and from this, the optimal control can be derived. In general, the HJB equation must be solved numerically. This solution has a problem of dimensionality as its complexity is exponential in the number of variables (state space). Some papers addressing the problem of dimensionality in connection to the HJB equation, and how to overcome it, are Garcke and Kr\"oner \cite{GK} and Kalise and Kunisch \cite{KK}. 
\end{remark}






From Peng \cite{Peng} (or {\O}ksendal and Sulem \cite{OS-Levy}, Section 3.2, for the generalised jump case), we have sufficient and necessary maximum principles for the solution of the stochastic control problem \eqref{eq: control_final}. 

For our PReP optimal stochastic control problem, the Hamiltonian is

\[
\begin{array}{llll}
\mathcal{H}(t, \bold{x}, u, \bold{p}, \bold{q}) &:=& w_1 x_2 + w_2 u^2 + b(\bold{x}, u)\bold{p} + \sigma{\bold{x}}\bold{q} \\[\medskipamount]
&=& w_1 x_2 + w_2 u^2 + (\bold{K} + f(\bold{x}) + A\bold{x} + B\bold{x}u)\bold{p} + g(\bold{x})\bold{q}.
\end{array}
\]

The adjoint processes $\bold{p}(t) = (p_1(t), \ldots, p_5(t))$, $\bold{q}(t) = (q_1(t), \ldots, q_5(t))$, $t \in [0,T]$ are given as solutions of the following system of BSDEs:

\begin{equation}
\label{eq: BSDE_adjoint}
\begin{array}{lll}
d\bold{p}(t) &=& -\nabla_{\bold{x}} \mathcal{H} (t)dt + \bold{q}(t)dB(t) \\[\medskipamount]
d \bold{p}(T) &=& 0
\end{array}
\end{equation}
where the terminal condition follows because there is no terminal time part in the performance function $J_{x_0}(u)$. That is,

\begin{equation}
\label{eq: adjoint}
\begin{array}{lllll}
dp_1(t) &=& [-\{{-}\beta (x_2 + \eta_C x_3 + \eta_A x_4) - (\mu + {\psi}) - u\}p_1(t) - \beta (x_2 + \eta_C x_3 + \eta_A x_4)p_2(t) \\[\smallskipamount]
&&- (u+ {\psi}) p_5(t) +\sigma(x_2+ \eta_C x_3 + \eta_A x_4)q_1(t) {-} \sigma(x_2+ \eta_C x_3 + \eta_A x_4)q_2(t)]dt \\[\smallskipamount]
&&+ q_1(t)dB(t) \\[\medskipamount]
dp_2(t) &=& [-w_1 {+} \beta x_1 p_1(t) - (\beta x_1 - \xi_3)p_2(t) - \phi p_3(t) - \rho p_4(t) + \sigma q_1(t)x_1 -\sigma q_2(t)x_1]dt \\[\medskipamount] 
&&+ q_2(t)dB(t) \\[\medskipamount]
dp_3(t) &=& [\beta x_1 \eta_C p_1(t)- \beta \eta_C x_1 p_2(t) - w p_2(t) +{\xi_2 }p_3(t) + \sigma \eta_C x_1 q_1(t) - \sigma \eta_C x_1 q_2(t)]dt \\[\smallskipamount]
&& + q_3(t) dB(t) \\[\medskipamount]
dp_4(t) &=& -\beta x_1 \eta_A p_1(t)- \beta \eta_A x_1 p_2(t) - \alpha p_2(t) +\xi_1 p_4(t) + \sigma \eta_A x_1 q_1(t) - \sigma \eta_A x_1 q_2(t)]dt \\[\smallskipamount]
&&+ q_4(t)dB(t) \\[\medskipamount]
dp_5(t) &=& [-\theta p_1(t) + \xi_4 p_5(t)]dt + q_5(t)dB(t) \\[\medskipamount]
p_i(T) &=& 0, \mbox{ } t=1, \ldots, 5.
\end{array}
\end{equation}




This system of BSDEs is linear. Hence, from Theorem 1.7 in {\O}ksendal and Sulem \cite{OksendalSulem}, we know that there exists a unique solution $\{(\bold{p}(t), \bold{q}(t))\}_{t \in [0,T]}$ to the system of adjoint BSDEs. Furthermore, both Theorem 1.7 in {\O}ksendal and Sulem from \cite{OksendalSulem} and Proposition 1.3. in El Karoui et al. from \cite{ElKarouiEtAl} provide an explicit solution to the BSDE. 

To derive the optimal control, we use the first order condition of the maximum principle:

\[
\begin{array}{lllll}
\frac{\partial \mathcal{H}(t, \bold{x}, u, \bold{p}, \bold{q})}{\partial u} &=&0, \mbox{ i.e.,} \\[\medskipamount]
2w_2u (t)+ B \bold{X}(t)\bold{p}(t) &=&0
\end{array}
\]
\noindent where the last equation follows from the definition of $B$. Solving this equation with respect to $u(t)$, we find 

\[
u(t) = \frac{S(t)(p_1(t) -p_5(t))}{2w_2}
\]
\noindent where $p_1(t), p_5(t)$ must be found by solving the adjoint BSDE system \eqref{eq: BSDE_adjoint}.

Note that this candidate optimal control is not necessarily in $\mathcal{A}$, since there is no guarantee that $u(t) = \frac{S(t)(p_1(t) -p_5(t))}{2w_2} \in [0,1]$ for almost all $t \in [0,T]$, $P$-a.s. However, if we instead consider

\begin{equation}
\label{eq: optimal_control_no_constraint}
\hat{u}(t) := \min \left\{ \max\left\{ 0,  \frac{S(t)(p_1(t) -p_5(t))}{2w_2}\right\},  1\right\},
\end{equation}
\noindent we can check that this candidate optimal stochastic control satisfies all the conditions of the sufficient maximum principle of Section 3.2 in {\O}ksendal and Sulem \cite{OS-Levy}. Hence, $\{\hat{u}(t)\}_{t \in [0,T]}$ is an optimal stochastic control. We summarise this result in a theorem:

\begin{theorem}
\label{thm: opt_control}
The stochastic optimal control, $\hat{u}(t)$, $t \geq 0$, corresponding to the PReP problem \eqref{eq: prep_problem} with SDE dynamics \eqref{eq: SDE_system} is given by:

\[
\hat{u}(t) := \min \left\{ \max\left\{ 0,  \frac{S(t)(p_1(t) -p_5(t))}{2w_2}\right\},  1\right\} \mbox{ for all } t \in [0,T].
\]
\end{theorem}

\dproof
See the previous derivation.
\fproof

At each time $t \in [0,T]$, the fraction $u(t)$ tells us how many percent of the susceptible individuals should be given PReP based on the current level of information. In Section \ref{sec: numerical}, we will illustrate this result numerically.

\section{Generalised Lagrange multiplier methods for stochastic optimal control}
\label{sec: Lagrange}

This section generalises the results of Section 2 in Dahl and Stokkereit \cite{DahlStokkereit}, and the framework is the same as in this paper, but adopted for our problem of stochastic control for the PReP vaccine.

We derive generalised Lagrange multiplier methods which can be combined with stochastic optimal control methods to solve the PReP stochastic control problem with either Type $I$  given with eq.(\ref{eq: type_i}) or Type $II$  given with eq.(\ref{eq: type_ii})  budget constraints. For these theoretical results, we consider the more general framework of stochastic jump processes. Hence, we consider a state process which may involve jumps, and a performance function with both an integral and a terminal time term. In Section \ref{sec: application}, we then apply the general results to the special case of PReP stochastic optimal control.

Consider the same framework as in Section \ref{sec: control_problem}, but in addition, let $\int_{\mb{R}} z \tilde{N}(dt,dz)$ a pure jump process independent of $B(t)$. Let $f : \mb{R}_+ \times \mb{R} \times \mb{R} \rightarrow \mb{R}$ and $g: \mb{R} \rightarrow \mb{R}$ be given, continuous functions. Let $\{\mathcal{F}_t\}_{t \in [0,T]}$ be the filtration generated by the Brownian motion and the pure jump process. We consider the stochastic optimal control problem which comes in two versions: Type $I$ and Type $II$ (see Section \ref{sec: control_problem}).
\begin{equation}
 \label{eq: OpprProblem}
\begin{array}{lll}
 \inf_{u \in \mathcal{A}} E^x[\int_0 ^T f(t,X(t),u(t)) dt + g(X(T))] \\[\smallskipamount]
\mbox{subject to} \\[\smallskipamount]
dX(t) = b(t,X(t),u(t))dt + \sigma(t,X(t),u(t))dB(t) \\[\smallskipamount]
\hspace{1.35cm} + \int_{\mb{R}} \gamma(t,X(t^-),u(t^-),z) \tilde{N}(dt,dz) \\[\smallskipamount]
(I) \mbox{ } E^x\left[\int_0^T M(t, X(t), u(t))dt\right] = 0 \mbox{ or } (II) \mbox{ } \int_0^T M(t, X(t), u(t))dt = 0 \mbox{ a.s.},
\end{array}
\end{equation}
\noindent where $M : \mb{R} \rightarrow \mb{R}$ is some given continuous function, $\mathcal{U} \subset \mb{R}$ is a given set, $b: \mb{R}_+ \times \mb{R} \times \mathcal{U} \rightarrow \mb{R}$, $\sigma: \mb{R}_+ \times \mb{R} \times \mathcal{U} \rightarrow \mb{R}$ and $\gamma: \mb{R}_+ \times \mb{R} \times \mathcal{U} \times \mb{R} \rightarrow \mb{R}$ . Here, $E^x[\cdot]$ denotes the expectation given that the state process $(\bold{X}(t))_{t \in [0,T]}$ starts in $x$, i.e. $\bold{X}(0)=x$.

In problem~(\ref{eq: OpprProblem}), the stochastic process $u(t) = u(t,\omega): \mb{R}_+ \times \Omega \rightarrow \mathcal{U}$ is our control process. We say that this control process $u(t)$ is admissible, and write $u \in \mathcal{A}$ if the dynamics of $X$ (i.e., the SDE in problem~(\ref{eq: OpprProblem})) has a unique, strong solution for all $x \in \mb{R}$, and
\[
E^x\left[\int_0 ^T f(t,X(t),u(t)) dt + g(X(T))\right] < \infty.
\]

\begin{remark}
The difference between problem \eqref{eq: OpprProblem} and the problem in Dahl and Stokkereit \cite{DahlStokkereit},  is that in \cite{DahlStokkereit}, the constraints are of the form 

\[
E^x\left[ M(X^u(T))\right] = 0 \mbox{ or } \int_0^T M(X^u(T))= 0 \mbox{ a.s.}.
\]
Hence, the current framework generalises that of Dahl and Stokkereit \cite{DahlStokkereit}.

\end{remark}

As seen in Section \ref{sec: unconstrained}, the stochastic maximum principle for jump diffusions by Framstad et al.~\cite{FOS} (see also Tang and Li~\cite{TangLi} and {\O}ksendal \cite{Oksendal}), can be used to find the optimal control of problem~(\ref{eq: OpprProblem}) without the constraints of Type $I$ or $II$. However, if we add a constraint such as $I$ or $II$ to the problem, such as in (\ref{eq: OpprProblem}), the stochastic maximum principle cannot be used directly. In this section, we show how the constrained stochastic optimal control problems can be solved by combining a generalised Lagrange duality method and the stochastic maximum principle.

For notational simplicity, problem~(\ref{eq: OpprProblem}) is assumed to be in one dimension. However, the results of this section also apply to multi-dimensional stochastic optimal control problems. The results generalise in a straight-forward manner (essentially just some extra notation). Also, note that even though our PReP control problem is without jump,  we include jumps in the framework for (\ref{eq: OpprProblem}). Since the jump framework is more general, the PReP problem is just a special case.

\subsection{Type I constraint}
\label{subsec: Expected}

Consider problem~(\ref{eq: OpprProblem}) with a Type $I$ constraint, i.e.:
\[
\begin{array}{ll}
\phi(x) \mbox{ }:= \mbox{ } \inf_{u \in \mathcal{A}} \mbox{ } E^x[\int_0 ^T f(t,X(t),u(t)) dt + g(X(T))] \hspace{2.5cm}\\[\smallskipamount]
\mbox{subject to}
\end{array}
\]
\begin{equation}
\begin{array}{rllll}
 \label{eq: Problem1}
dX(t) &=&b(t,X(t),u(t))dt + \sigma(t,X(t),u(t))dB(t) \\[\smallskipamount]
&&+ \int_{\mb{R}} \gamma(t,X(t^-),u(t^-),z) \tilde{N}(dt,dz) \\[\smallskipamount]
E^x[\int_0^T M(t, X(t), u(t))dt]  &=& 0.
\end{array}
\end{equation}




This problem can be solved using the standard Lagrange multiplier method, and then applying some method of stochastic control, for instance the stochastic maximum principle. Hence, let $\lambda \in \mb{R}$ be a Lagrange multiplier. Then, we introduce the unconstrained stochastic control problem
\begin{equation}
 \label{eq: Problem1Lagrange}
\begin{array}{llllll}
\phi_{\lambda}(x) \mbox{ } := \mbox{ } \inf_{u \in \mathcal{A}} \mbox{ } E^x[\int_0 ^T f(t,X(t),u(t)) dt + g(X(T)) + \lambda \int_0^T M(t, X(t), u(t))dt] \\[\smallskipamount]
\mbox{subject to} \\[\smallskipamount]
\hspace{1.6cm}dX(t) = b(t,X(t),u(t))dt + \sigma(t,X(t),u(t))dB(t) \\[\smallskipamount]
\hspace{3cm}+ \int_{\mb{R}} \gamma(t,X(t^-),u(t^-),z) \tilde{N}(dt,dz).
\end{array}
\end{equation}
This solution strategy is explored in Section 11.3 in {\O}ksendal~\cite{Oksendal} for the no-jump case. However, the proof of this theorem generalises in a straight-forward manner to the jump case. Therefore, we have the following theorem.




\begin{theorem}
\label{thm: 1}
(Type I, equality constraint)

\medskip

\noindent Suppose that we for all $\lambda \in \mb{R}$ one can find $\phi_{\lambda}(y)$ and $u_{\lambda}^*$ solving the unconstrained stochastic control problem~(\ref{eq: Problem1Lagrange}). Moreover, suppose there exists $\lambda_0 \in \mb{R}$ such that
\[
 E^x\left[\int_0^T M(t, X_{u_{\lambda_0}^*}(t), u_{\lambda_0}^*(t))dt\right] =0.
\]
Then, $\phi(x) := \phi_{\lambda_0}(x)$ and $u^* := u^*_{\lambda_0}$ solves the constrained stochastic control problem~(\ref{eq: Problem1}).
\end{theorem}

\dproof
 See {\O}ksendal~\cite{Oksendal}, Theorem 11.3.1. The proof is also similar to the proof of the following Theorem~\ref{thm: 2}.
\fproof




From Theorem~\ref{thm: 1}, in order to solve problem~(\ref{eq: OpprProblem}), it is sufficient to solve problem~(\ref{eq: Problem1Lagrange}), and then determine a Lagrange multiplier $\lambda_0$ which satisfies the Type I constraint for this optimal control. Note that problem~(\ref{eq: Problem1Lagrange}) can be solved using the stochastic maximum principle.




\begin{remark}
 \label{remark: future work-SDEwthjumps}
For the sake of generality of the theory, the set up (\ref{eq: Problem1Lagrange}) includes the possibility of jumps, even-though the system of SDEs which describes the spread of HIV with PReP treatment is without jumps. An idea for future work is to introduce jumps in the SDE model (\ref{eq: SDE_system}). Considering the PReP model with jumps may be more realistic than the current no-jump case, due to the possibility of disasters and crises, economical or natural, which may influence the number of infected and susceptible individuals. For this kind of model, a complete analysis of the system of SDEs should be obtained, as done in the no-jump case Djordevi{\'c} and Silva  in \cite{DjordevicSilva}.
\end{remark}

\subsection{Type II constraint}
\label{subsec: StokLagrLik}

Now, consider problem~(\ref{eq: OpprProblem}) with a type $II$ constraint:

\[
\begin{array}{lll}
\phi(x) \mbox{ } := \mbox{ } \inf_{u \in \mathcal{A}} \mbox{ } E^x[\int_0 ^T f(t,X(t),u(t)) dt + g(X(T))] \hspace{2.5cm} \\[\smallskipamount]
\mbox{subject to} \\[\smallskipamount]
\end{array}
\]
\begin{equation}
 \label{eq: Problem2}
\begin{array}{lrllll}
&dX(t) &=&b(t,X(t),u(t))dt + \sigma(t,X(t),u(t))dB(t) \\[\smallskipamount]
&&&+ \int_{\mb{R}} \gamma(t,X(t^-),u(t^-),z) \tilde{N}(dt,dz) \\[\smallskipamount]
&\int_0^T M(t, X(t), u(t))dt &=& 0 \mbox{ a.s.}
\end{array}
\end{equation}

\noindent where, as before, $M : \mb{R} \rightarrow \mb{R}$ is a given, continuous function. For notational simplicity, let us define the performance function
%




\[
J^u(x) := E^x\left[\int_0 ^T f(t,X(t),u(t)) dt + g(X(T))\right].
\]
We would like to use the Lagrange multiplier concept to solve problem~(\ref{eq: Problem2}) by solving an unconstrained stochastic control problem. However, since we have an almost sure constraint, it is not sufficient to introduce a single scalar Lagrange multiplier $\lambda \in \mb{R}$. The Lagrange multiplier must be stochastic in order to handle the stochastic constraint $\int_0^T M(t, X(t), u(t))dt = 0 \mbox{ a.s.}$ Hence, we introduce an $\mathcal{F}_T$-measurable stochastic Lagrange multiplier $\mu : \Omega \rightarrow \mb{R}$ (which we will also call a \emph{stochastic multiplier}). Note that $\mu$ must be $\mathcal{F}_T$-measurable, since $\int_0^T M(t, X(t), u(t))dt$ is $\mathcal{F}_T$-measurable.

Assume that the stochastic multiplier $\mu$ satisfies $E[\mu] < \infty$. Moreover, assume that $E^x[\int_0^T M(t, X(t), u(t))dt < \infty$ for all $u \in \mathcal{A}$. We introduce a new stochastic control problem
\begin{samepage}
\[
\begin{array}{ll}
\phi_{\mu}(x) \mbox{ } := \mbox{ }  \inf_{u \in \mathcal{A}} \mbox{ } E^x[\int_0 ^T f(t,X(t),u(t)) dt + g(X(T)) + \mu \int_0^T M(t, X(t), u(t))dt] \hspace{1cm} \\[\medskipamount]
\mbox{subject to}
\end{array}
\]
\begin{equation}
 \label{eq: LagrangeProblem2}
\begin{array}{lrllll}
dX(t) &=&b(t,X(t),u(t))dt + \sigma(t,X(t),u(t))dB(t) \\[\smallskipamount]
&&+ \int_{\mb{R}} \gamma(t,X(t^-),u(t^-),z) \tilde{N}(dt,dz),
\end{array}
\end{equation}
\end{samepage}

\noindent and define
\[
J^u_{\mu}(x) := E^x\left[\int_0 ^T f(t,X(t),u(t)) dt + g(X(T)) + \mu \int_0^T M(t, X(t), u(t))dt\right].
\]

We also define the set of stochastic multipliers by

\[
 \Lambda := \{ \mu : \Omega \rightarrow \mb{R}\mbox{ } | \mbox{ } \mu \mbox{ is } \mathcal{F}_T\mbox{-measurable and } E[\mu] < \infty\}.
\]

The following Theorem~\ref{thm: 2} states that if there exists a solution to the unconstrained problem~\eqref{eq: LagrangeProblem2} with a stochastic multiplier which ensures that the constraint $\int_0^T M(t, X(t), u(t))dt = 0$ a.s. is satisfied, then we have a corresponding solution to our original problem~\eqref{eq: Problem2}.




\begin{theorem}
 \label{thm: 2}
 (Type II, equality constraint)
Suppose that we for all $\mu \in \Lambda$ can find $\phi_{\mu}(x)$ and $u_{\mu}^*$ solving the unconstrained stochastic control problem~(\ref{eq: LagrangeProblem2}). Moreover, suppose there exists $\mu_0 \in \Lambda$ such that
\[
\int_0^T M(t, X_{u_{\mu_0}^*}(t), u_{\mu_0}^*(t)dt) =0 \mbox{ a.s.}
\]
Then, $\phi(x) := \phi_{\mu_0}(x)$ and $u^* := u^*_{\mu_0}$ solves the constrained stochastic control problem~(\ref{eq: Problem2}).

\end{theorem}

\dproof
Let $\mu$ be $\mathcal{F}_T$-measurable. Then,

\[
\begin{array}{lllll}
 E^x[\int_0^T f(t,u_{\mu}^*, X_{u_{\mu}^*})dt + g(X_{u_{\mu}^*}(T)) + \mu \int_0^T M(t,X_{u_{\mu}^*}(t), u_{\mu}^*(t))dt] = J_{\mu}^{u_{\mu}^*}(x) \\[\medskipamount]
\leq J^{u}_{\mu}(x) = E^x[\int_0^T f(t,u, X_{u})dt + g(X_{u}(T)) + \mu \int_0^T M(t, X_{u}(t), u(t))dt]
\end{array}
\]
\noindent where the first equality uses the definition of $J^u_\mu$, the inequality uses the definition of $u_{\mu}^*$ and the final equality uses the definition of $J^u_\mu$.

In particular, if $\mu = \mu_0$ a.s. and $u$ is feasible in the constrained control problem~(\ref{eq: Problem2}), then

\begin{equation}
\label{eq: books}
\int_0^T M(t, X_{u_{\mu_0}^*}(t), u_{\mu_0}^*(t))dt =0= \int_0^T M(t, X_{u}(t), u(t))dt \mbox{ a.s.}
\end{equation}
\noindent from the definition of $\mu_0$ and the assumption that $u$ is feasible in problem~(\ref{eq: Problem2}).

Hence,
\[
\begin{array}{lll}
 E^x[\int_0^T f(t,u_{\mu_0}^*, X_{u_{\mu_0}^*})dt + g(X_{u_{\mu_0}^*}(T)) + \mu_0 \int_0^T M(t, X_{u_{\mu_0}^*}(t), u_{\mu_0}^*(t))dt] \\[\smallskipamount]=J^{u^*_{\mu_0}}_{\mu_0}(x)
\leq J^{u}_{\mu_0}(x) =  E^x[\int_0^T f(t,u, X_{u})dt + g(X_{u}(T)) + \mu_0 \int_0^T M(t, X_{u}(t), u(t))dt]
\end{array}
\]

As  $\int_0^T M(t, X_{u_{\mu_0}^*}(t), u_{\mu_0}^*(t))dt =0= \int_0^T M(t, X_{u}(t), u(t))dt \mbox{ a.s.}$ from equation~(\ref{eq: books}), so

\[
J^{u^*_{\mu_0}}(x) = J^{u^*_{\mu_0}}_{\mu_0}(x) \leq J^{u}_{\mu_0}(x) = J^{u}(x)
\]
\noindent for all stochastic controls $u$ feasible in the constrained problem~(\ref{eq: Problem2}). Note that $u^*_{\mu_0}$ is feasible in problem~(\ref{eq: Problem2}), therefore it is an optimal control for this problem.
\fproof

\medskip

Note that problem~(\ref{eq: LagrangeProblem2}) is a stochastic optimal control problem of the form in {\O}ksendal and Sulem~\cite{OS-Levy}, with $f_{\mu}(\cdot)=f(\cdot) + \mu M(\cdot)$ and $g_{\mu}(\cdot)=g(\cdot)$. Therefore, we may use some known methods of stochastic control, for example the stochastic maximum principle, to solve the problem. Note that it is irrelevant for this solution strategy whether the unconstrained stochastic control problem coming from the stochastic Lagrange multiplier method is solved using the maximum principle, or some other method of stochastic control. If it is more suitable for the problem at hand, the dynamic programming/Hamilton-Jacobi-Bellman approach to stochastic control of jump diffusions can also be used, see {\O}ksendal and Sulem~\cite{OS-Levy} Theorem 3.1. For the dynamic programming approach, the problem must have a Markovian structure.




Theorem \ref{thm: 1} and Theorem \ref{thm: 2} both consider equality constraints, however, in the PReP stochastic optimal control problem, constraints are defined with the inequalities  of the form:

\[
\begin{array}{llll}
\mbox{Type $I$: } E[\int_0^T S(t)u(t) c(t) dt] - \mathcal{V}&\leq& 0, \\[\medskipamount]

\mbox{Type $II$: } \int_0^T S(t)u(t) c(t) dt -\mathcal{V}&\leq& 0 \mbox{ } P\mbox{-a.s.} 
\end{array}
\]

\noindent The generalised version of these kinds of constraints are:

\[
\begin{array}{llll}
\mbox{Type $I$: } E[\int_0^T M(t, X(t), u(t))) dt]&\leq& 0, \\[\medskipamount]

\mbox{Type $II$: } \int_0^T M(t, X(t), u(t))dt &\leq& 0 \mbox{ } P\mbox{-a.s.}\ .
\end{array}
\]

\noindent However, Theorem \ref{thm: 1} and Theorem \ref{thm: 2} can be generalised to the inequality case in a straight-forward manner, by simply adding that the Lagrange multipliers have to be non-negative ($P$-a.s. in Theorem \ref{thm: 2}). The inequality constraint versions of the two theorems are given in the sequel.

\begin{theorem}
\label{thm: 1_ineq}
(Type I, inequality constraint)
\noindent Consider the stochastic optimal control problem with a Type I inequality constraint. Suppose that we for all $\lambda \in \mb{R}$ can find $\phi_{\lambda}(y)$ and $u_{\lambda}^*$ solving the unconstrained stochastic control problem~(\ref{eq: Problem1Lagrange}). Moreover, suppose there exists $\lambda_0 \geq 0$ such that
\[
 E^x\left[\int_0^T M(t, X_{u_{\lambda_0}^*}(t), u_{\lambda_0}^*(t))dt\right] = 0.
\]
Then, $\phi(x) := \phi_{\lambda_0}(x)$ and $u^* := u^*_{\lambda_0}$ solves the constrained stochastic control problem~(\ref{eq: Problem1}).
\end{theorem}

\dproof
This is a straight forward generalization of Theorem \ref{thm: 1}, and therefore, we omit writing it out again.
\fproof

\begin{theorem}
 \label{thm: 2_ineq}
 (Type II, inequality constraint)
 \noindent Consider the stochastic optimal control problem with a Type II inequality constraint. Suppose that we for all $\mu \in \Lambda$ can find $\phi_{\mu}(x)$ and $u_{\mu}^*$ solving the unconstrained stochastic control problem~(\ref{eq: LagrangeProblem2}). Moreover, suppose there exists $\mu_0 \in \Lambda$ with $\mu_0 \geq 0$ $P$-a.e., such that
\[
\int_0^T M(t, X_{u_{\mu_0}^*}(t), u_{\mu_0}^*(t)dt)  = 0 \mbox{ a.s.}
\]
Then, $\phi(x) := \phi_{\mu_0}(x)$ and $u^* := u^*_{\mu_0}$ solves the constrained stochastic control problem~(\ref{eq: Problem2}).

\end{theorem}


\dproof

Let $\mu$ be $\mathcal{F}_T$-measurable. Then,

\[
\begin{array}{lllll}
 E^x[\int_0^T f(t,u_{\mu}^*, X_{u_{\mu}^*})dt + g(X_{u_{\mu}^*}(T)) + \mu \int_0^T M(t,X_{u_{\mu}^*}(t), u_{\mu}^*(t))dt] = J_{\mu}^{u_{\mu}^*}(x) \\[\medskipamount]
\leq J^{u}_{\mu}(x) = E^x[\int_0^T f(t,u, X_{u})dt + g(X_{u}(T)) + \mu \int_0^T M(t, X_{u}(t), u(t))dt]
\end{array}
\]
\noindent where the first equality uses the definition of $J^u_\mu$, the next one uses the definition of $u_{\mu}^*$ and the final equality uses the definition of $J^u_\mu$.

In particular, if $\mu = \mu_0$ a.s. and $u$ is feasible in the constrained control problem, then

\begin{equation}
\label{eq: boks}
\int_0^T M(t, X_{u_{\mu_0}^*}(t), u_{\mu_0}^*(t))dt = 0 \mbox { and } \int_0^T M(t, X_{u}(t), u(t))dt \leq 0 \mbox{ a.s.}
\end{equation}
\noindent from the definition of $\mu_0$ and the assumption that $u$ is feasible in problem.

Hence,
\[
\begin{array}{lll}
J^{u^*_{\mu_0}}_{\mu_0}(x) \leq J^{u}_{\mu_0}(x) 
\end{array}
\]

But as $\int_0^T M(t, X_{u_{\mu_0}^*}(t), u_{\mu_0}^*(t))dt =0$ and  $\int_0^T M(t, X_{u}(t), u(t))dt  \leq 0 \mbox{ a.s.}$ (eq. (\ref{eq: boks})), and since $\mu_0 \geq 0$ $P$-a.s., it follows that

 
 

\[
J^{u^*_{\mu_0}}(x) = J^{u^*_{\mu_0}}_{\mu_0}(x) \leq J^{u}_{\mu_0}(x){ \leq} J^{u}(x)
\]



\noindent for all stochastic controls $u$ feasible in the constrained problem.

Note that $u^*_{\mu_0}$ is feasible in the constrained control problem and therefore it is an optimal control for this problem.
\fproof

\section{Optimal control of PReP with budget constraint}
\label{sec: application}

In this section, Theorem \ref{thm: 1} and \ref{thm: 2}  will be applied in order to solve the PReP stochastic optimal control problem with constraints of 
Type $I$ and Type $II$, respectively. The framework for the PReP-problem is slightly simpler than the framework in Section \ref{sec: Lagrange} because we don't have any jump terms and the terminal time term of the performance function is zero. The constraints which are consider are:

\[
\begin{array}{llll}
\mbox{Type I: } E[\int_0^T S(t)u(t) c(t) dt] &\leq& \mathcal{V},\\[\medskipamount]

\mbox{Type II: } \int_0^T S(t)u(t) c(t) dt &\leq& \mathcal{V} \mbox{ } a.s. \,
\end{array}
\]
\noindent where $c(t)$ is some given cost function for PReP. That is, $c(t, \omega)$, $t \in [0,T], \omega \in \Omega$ is the cost of a single individual being treated with PReP. In other words, the Type $I$ constraint states that the expected total cost of PReP treatment over the whole time period of interest should not exceed $\mathcal{V}$. The Type $II$ constraint states that the total cost of PReP treatment over the whole time period of interest should not exceed $\mathcal{V}$ almost surely. As previously mentioned, the Type $II$ constraint is stricter than the Type $I$ constraint. 

Note that both of these constraints are slightly different from the one considered in Silva and Torres \cite{SilvaTorres}. They consider a deterministic constraint of the form $S(t)u(t) \leq \mathcal{V}$ a.s. In words, the total number of individuals treated with PReP should never exceed the pre-determined level $\mathcal{V}$. The reason we consider an integral constraint instead is that the constraint considered in Silva and Torres \cite{SilvaTorres} cannot be handled via the Lagrange techniques of Section \ref{sec: Lagrange}. Furthermore, we believe that the integral budget constraint is as realistic as constraining the number of treated individuals at any time.



\subsection{Type I constraint}
\label{sec: type_i}

Let $\lambda \geq 0$ be deterministic. The unconstrained Lagrange version of the stochastic control problem is

\begin{equation}
\label{eq: control_final}
\begin{array}{lll}
\min_{u \in \mc{A}} E^{x_0} [\int_0^T w_1 I(t) + w_2 u^2(t)dt + \lambda \int_0^T S(t)u(t)c(t)dt] \\[\medskipamount]
\mbox{s.t.} \\[\medskipamount]
d\bold{X}(t)= b(\bold{X}(t), u(t))dt + \sigma(\bold{X}(t))dB(t), \mbox{ } t \in [0,T].
\end{array}
\end{equation}

\noindent The Hamiltonian is
\[
\begin{array}{llll}
\mathcal{H}(t, \bold{x}, u, \bold{p}, \bold{q}, \lambda) 
&=& w_1 x_2 + w_2 u^2 + \lambda x_1 u c+ (\bold{K} + f(\bold{x}) + A\bold{x} + B\bold{x}u)\bold{p} + g(\bold{x})\bold{q}.
\end{array}
\]

The adjoint processes $\bold{p}(t) = (p_1(t), \ldots, p_5(t))$, $\bold{q}(t) = (q_1(t), \ldots, q_5(t))$, $t \in [0,T]$ are given as solutions of the following system of BSDEs:

\begin{equation}
\label{eq: BSDE_adjoint}
\begin{array}{lll}
d\bold{p}(t) &=& -\frac{\partial \mathcal{H}}{\partial \bold{x}}(t)dt + \bold{q}(t)dB(t) \\[\medskipamount]
d \bold{p}(T) &=& 0
\end{array}
\end{equation}
where the terminal condition follows because there is no terminal time part in the performance function. That is,

\[
\begin{array}{lllll}
dp_1(t) &=& [-\lambda u c-\{-\beta (x_2 + \eta_C x_3 + \eta_A x_4) - (\mu + {\psi}) - u\}p_1(t) - \beta (x_2 + \eta_C x_3 + \eta_A x_4)p_2(t) \\[\smallskipamount]
&&- (u+  {\psi}) p_5(t) + \sigma(x_2+ \eta_C x_3 + \eta_A x_4)q_1(t) -  \sigma(x_2+ \eta_C x_3 + \eta_A x_4)q_2(t)]dt \\[\smallskipamount]
&&+ q_1(t)dB(t) \\[\medskipamount]
dp_2(t) &=& [-w_1  {+}\beta x_1 p_1(t) - (\beta x_1 - \xi_3)p_2(t) - \phi p_3(t) - \rho p_4(t) + \sigma q_1(t) {x_1} -\sigma q_2(t)x_1]dt \\[\medskipamount] 
&&+ q_2(t)dB(t) \\[\medskipamount]
dp_3(t) &=& [\beta x_1 \eta_C p_1(t)- \beta \eta_C x_1 p_2(t) - w p_2(t) + {\xi_2} p_3(t) + \sigma \eta_C x_1 q_1(t) - \sigma \eta_C x_1 q_2(t)]dt \\[\smallskipamount]
&& + q_3(t) dB(t) \\[\medskipamount]
dp_4(t) &=& [-\beta x_1 \eta_A p_1(t)- \beta \eta_A x_1 p_2(t) - \alpha p_2(t) +\xi_1 p_4(t) + \sigma \eta_A x_1 q_1(t) - \sigma \eta_A x_1 q_2(t)]dt \\[\smallskipamount]
&&+ q_4(t)dB(t) \\[\medskipamount]
dp_5(t) &=& [-\theta p_1(t) + \xi_4 p_5(t)]dt + q_5(t)dB(t) \\[\medskipamount]
p_i(T) &=& 0, \mbox{ } t=1, \ldots, 5.
\end{array}
\]




This system of backward stochastic differential equations  (BSDEs) is linear. Hence, from Theorem 1.7 in {\O}ksendal and Sulem \cite{OksendalSulem}, we know that there exists a unique solution $\{(\bold{p}(t), \bold{q}(t))\}_{t \in [0,T]}$ to the system of adjoint BSDEs. This theorem also gives an explicit formula for the solution to the BSDE.

To derive the optimal control, we use the first order condition of the maximum principle:

\[
\begin{array}{lllll}
\frac{\partial \mc{H}(t, \bold{x}, u, \bold{p}, \bold{q})}{\partial u} &=&0, \mbox{ i.e.,} \\[\medskipamount]
2w_2u (t) + \lambda X_1(t)c(t)+ B \bold{X}(t)\bold{p}(t) &=&0,
\end{array}
\]
\noindent where the last equation follows from the definition of $B$. Solving this equation with respect to $u(t)$, we find 

\[
u(t) = \frac{S(t)(p_1(t) -p_5(t) - \lambda c(t))}{2w_2},
\]
\noindent where $p_1(t), p_5(t)$ must be found by solving the adjoint BSDE system \eqref{eq: BSDE_adjoint}.

Note that this candidate optimal control is not necessarily in $\mathcal{A}$, since there is no guarantee that $u(t) = \frac{S(t)(p_1(t) -p_5(t)-\lambda c(t))}{2w_2} \in [0,1]$ for almost all $t \in [0,T]$, $P$-a.s. However, if we instead consider

\begin{equation}
\label{eq: optimal_control_typei_constraint}
u^*(t) := \min \left\{ \max\{ 0,  \frac{S(t)(p_1(t) -p_5(t) - \lambda c(t))}{2w_2}\},  1\right\}.
\end{equation}
\noindent 
It could be easily  checked that this candidate optimal stochastic control satisfies all the conditions of the sufficient maximum principle of Section 3.2 in {\O}ksendal and Sulem \cite{OS-Levy}. Hence, by Theorem \ref{thm: 1_ineq} and the stochastic maximum principle, $\{u^*(t)\}_{t \in [0,T]}$ is an optimal stochastic control if there exists a $\lambda_0 \geq 0$ such that

\[
E\left[\int_0^T S_{u_{\lambda_0}^*}(t) u_{\lambda_0}^*(t) c(t)dt\right]  - \mathcal{V}= 0.
\]




\subsection{Type II constraint}
\label{sec: type_ii}

For a Type II inequality constraint, the problem formulation becomes identical to the one in Section \ref{sec: type_i}, except that the Lagrange multiplier $\lambda$ is stochastic. The derivation is also identical,

\begin{equation}
\label{eq: optimal_control_typei_constraint}
u^*(t) := \min \left\{ \max\{ 0,  \frac{S(t)(p_1(t) -p_5(t) - \lambda c(t))}{2w_2}\},  1\right\},
\end{equation}
\noindent we can check that this candidate optimal stochastic control satisfies all the conditions of the sufficient maximum principle. Hence, by Theorem \ref{thm: 2_ineq} and the stochastic maximum principle, $\{u^*(t)\}_{t \in [0,T]}$ is an optimal stochastic control if there exists a random variable $\Lambda_0$ such that $\Lambda_0 \geq 0$ $P$-a.e. and

\[
\int_0^T S_{u_{\Lambda_0}^*}(t) u_{\Lambda_0}^*(t) c(t)dt  - \mathcal{V}= 0 \mbox{ } P\mbox{-a.e}.
\]

\section{Numerical example: The unconstrained stochastic control of PReP problem}
\label{sec: numerical}

In this section, we present a numerical example to show a practical application of the previous unconstrained stochastic optimal control problem from Section \ref{sec: unconstrained}. We extend the numerical example from Section \ref{sec: prep_sde}, and use the parameter values in Table \ref{fig: table}. We generalise the method from Campos et al. \cite{Campos} to the stochastic case: In order to derive a numerical solution to the optimal control problem via Theorem \ref{thm: opt_control}, we need to iteratively solve the PReP controlled SDE dynamics \eqref{eq: SDE_system} and corresponding adjoint equation \eqref{eq: adjoint}. For each Monte Carlo path, this is done in as follows: 

While the absolute error of convergence is too large (in some sense, to be discussed later), repeat these 4 steps:

\begin{enumerate}
\item{Iterate the PReP controlled SDE \eqref{eq: SDE_system} via a forward stochastic first order Runge-Kutta method. To initialise the method, we use a guess for the control $u$ (in our case, $u=0$) and choose some initial conditions $(S(0), I(0), A(0), C(0), E(0))$ (see Table \ref{fig: table}). This results in an approximate solution to the PReP controlled SDE  \eqref{eq: SDE_system}, similar to that presented in Section \ref{sec: prep_sde}: $(\hat{S}, \hat{I}, \hat{A}, \hat{C}, \hat{E})$.}
\item{The approximate solution $(\hat{S}, \hat{I}, \hat{A}, \hat{C}, \hat{E})$ is then used as input to solve the adjoint BSDE \eqref{eq: adjoint}. The adjoint BSDE is solved via a backward stochastic first order Runge-Kutta method. This results in an approximate solution to the adjoint equation: $(\hat{p}_1, \hat{p}_2, \hat{p}_3, \hat{p}_4, \hat{p}_5)$.}
\item{Then, the control $u$ is updated by taking a convex combination of the previous iteration of the control and a new control computed from the formula in Theorem \ref{thm: opt_control}. That is,
\[
u_i = \lambda_1 u^{\mbox{old}}_{i} + \lambda_2 \hat{u_i}
\]
where $\hat{u}_i$ is from the formula in Theorem \ref{thm: opt_control}:
\[
\hat{u}_i := \min \left\{ \max\left\{ 0,  \frac{\hat{S}_i(\hat{p}_{1,i} -\hat{p}_{5,i})}{2w_2}\right\},  1\right\} \mbox{ for all } t \in [0,T]
\]
\noindent and $u^{\mbox{old}}_i$ is from the previous pass of the while-loop. Here, $\lambda_1, \lambda_2 \in [0,1]$ are convex coefficients weighting old (i.e. the previous iteration) vs. new knowledge (i.e. the newest iteration). Note that the choice of $\lambda$ bears some resemblance to the so-called learning rate in machine learning.
}
\item{The while loop is stopped when the difference between the new iteration and the previous one is sufficiently small for all the processes 
\[
(S, I, A, C, E, p_1, p_2, p_3, p_4, p_5, q_1, q_2, q_3, q_4, q_5)
\] 
\noindent in comparison to the absolute size of the respective processes.}
\end{enumerate}

Further details about the numerics will be the topic of a forthcoming paper. For the purpose of this paper, we include the numerical example to illustrate our theoretical results about stochastic optimal control of the PReP problem.

\begin{figure}
\begin{centering}
\includegraphics[width=0.75\textwidth]{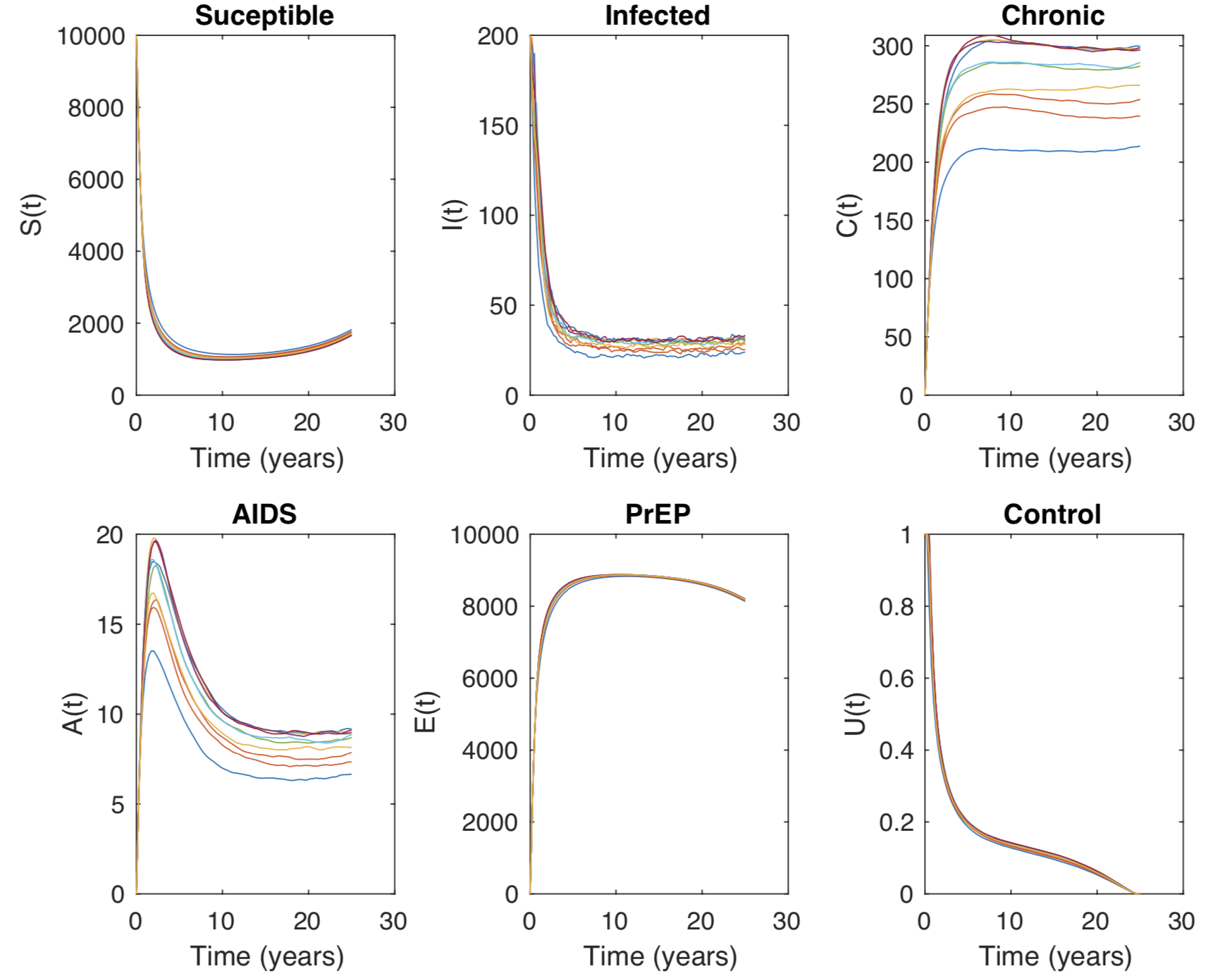}
\end{centering}
\caption{Plot of 10 paths of the solution of the stochastic optimal control PReP problem over 25 years with $\sigma =0.2/N$ and $N=10 200$.}
\label{fig: opt_control}
\end{figure}

In Figure \ref{fig: opt_control}, we have plotted 10 paths of the PReP process dynamics under the optimal control given by Theorem \ref{thm: opt_control} and found by the scheme above. The terminal time is $25$ years. In Figure \ref{fig: opt_control}, $\sigma = 0.2/N$. Like in Section \ref{sec: prep_sde}, we have initial conditions 

\[
S(0) = 10000,  I(0) = 200, C(0) = 0, A(0) = 0, E(0) = 0.
\]

We consider $N=10 200$, so $N=S(0) + I(0) + A(0) + C(0) + E(0)$. For the performance function, \eqref{eq: performance}, we have chosen the weights 
\[
      w_1 = 20, w_2 = 0.3N.
\]

This choice of weight implies that we weight avoiding infected individuals greatly in comparison to the cost of PReP.

\begin{figure}
\begin{centering}
\includegraphics[width=0.75\textwidth]{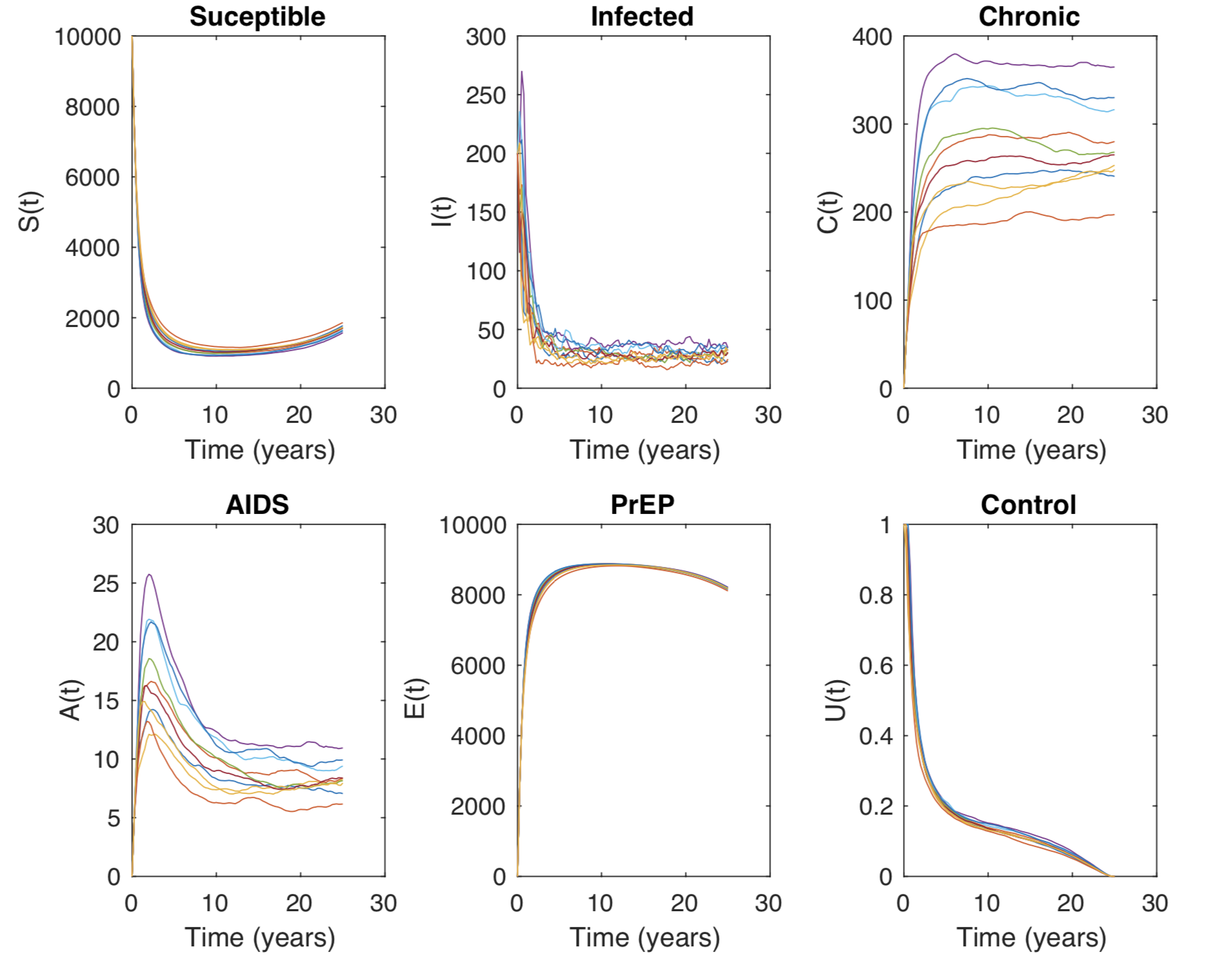}
\end{centering}
\caption{Plot of 10 paths of the solution of the stochastic optimal control PReP problem over 25 years with $\sigma =0.6/N$ and $N=10 200$.}
\label{fig: opt_control2}
\end{figure}

In Figure \ref{fig: opt_control2}, we have plotted 10 paths of the PReP process dynamics under the optimal control given by Theorem \ref{thm: opt_control} and found by the scheme above. In Figure \ref{fig: opt_control}, $\sigma = 0.6/N$ and we consider $N=10 200$. The only change from Figure \ref{fig: opt_control} is the value of $\sigma$. A larger $\sigma$, as in Figure \ref{fig: opt_control2} corresponds to a greater weighting of the noise terms in the model. The increase of noise is seen by comparing Figure \ref{fig: opt_control} and Figure \ref{fig: opt_control2}. The variance of the processes plotted in Figure \ref{fig: opt_control2} appear to be slightly larger than that in Figure $\ref{fig: opt_control}$.

Note that in both Figure \ref{fig: opt_control} and Figure \ref{fig: opt_control2}, the shape of the stochastic optimal control is essentially the same, and the variance is very low. According to the optimal control, one should start out with a very high PReP treatment rate (the maximum amount of PReP one allows for, in this case u=1). Then, the PReP treatment rate should rapidly decrease with time for approximately the first 3 years, before gradually being reduced for the remaining 22 years and ending up with no PReP treatment. The reason for this fading out of the PReP treatment is that the terminal time in the optimal control problem is set to be $T=25$. By recalling the performance function

\begin{equation}
\label{eq: performance}
J(u) := E\left[\int_0^T [w_1 I(t) + w_2 u^2(t)]dt\right],
\end{equation} 
\noindent we see that there is no weight added to what happens after the terminal time or at the actual terminal time. Because of this, it will be optimal to let the PReP treatment rate go towards zero as one approaches the terminal time.

\begin{figure}
\begin{centering}
\includegraphics[width=0.75\textwidth]{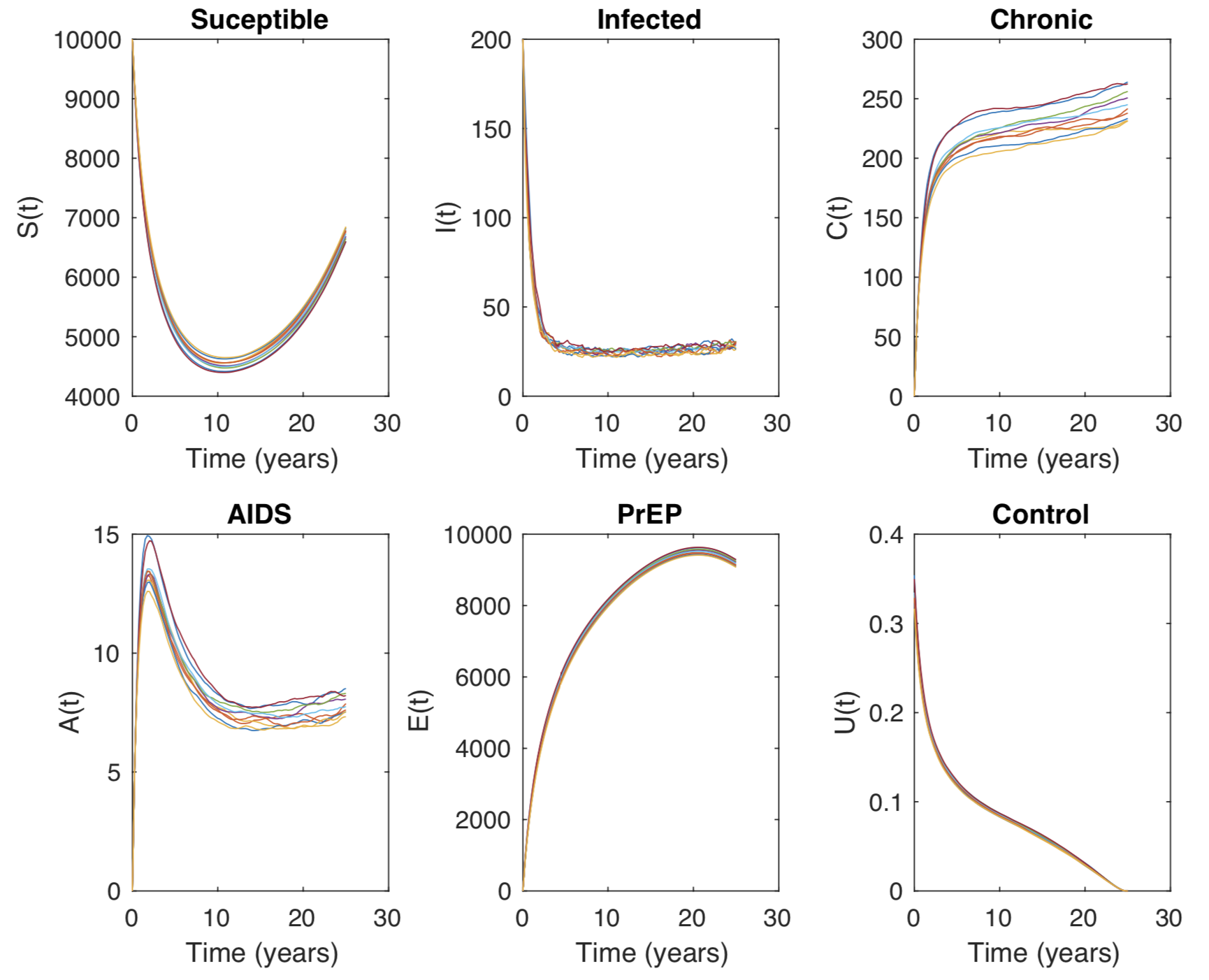}
\end{centering}
\caption{Plot of 10 paths of the solution of the stochastic optimal control PReP problem over 25 years with $\sigma = 0.2/N$ and $N=30 000$.}
\label{fig: opt_control3}
\end{figure}

In Figure \ref{fig: opt_control3}, we have again plotted 10 paths of the PReP model under the stochastic optimal control. The framework is as in Figure \ref{fig: opt_control}, except that we have chosen $N=30 000$. Hence, $N > S(0) + I(0) + A(0) + C(0) + E(0) = 10 200$. 

\begin{figure}
\begin{centering}
\includegraphics[width=0.75\textwidth]{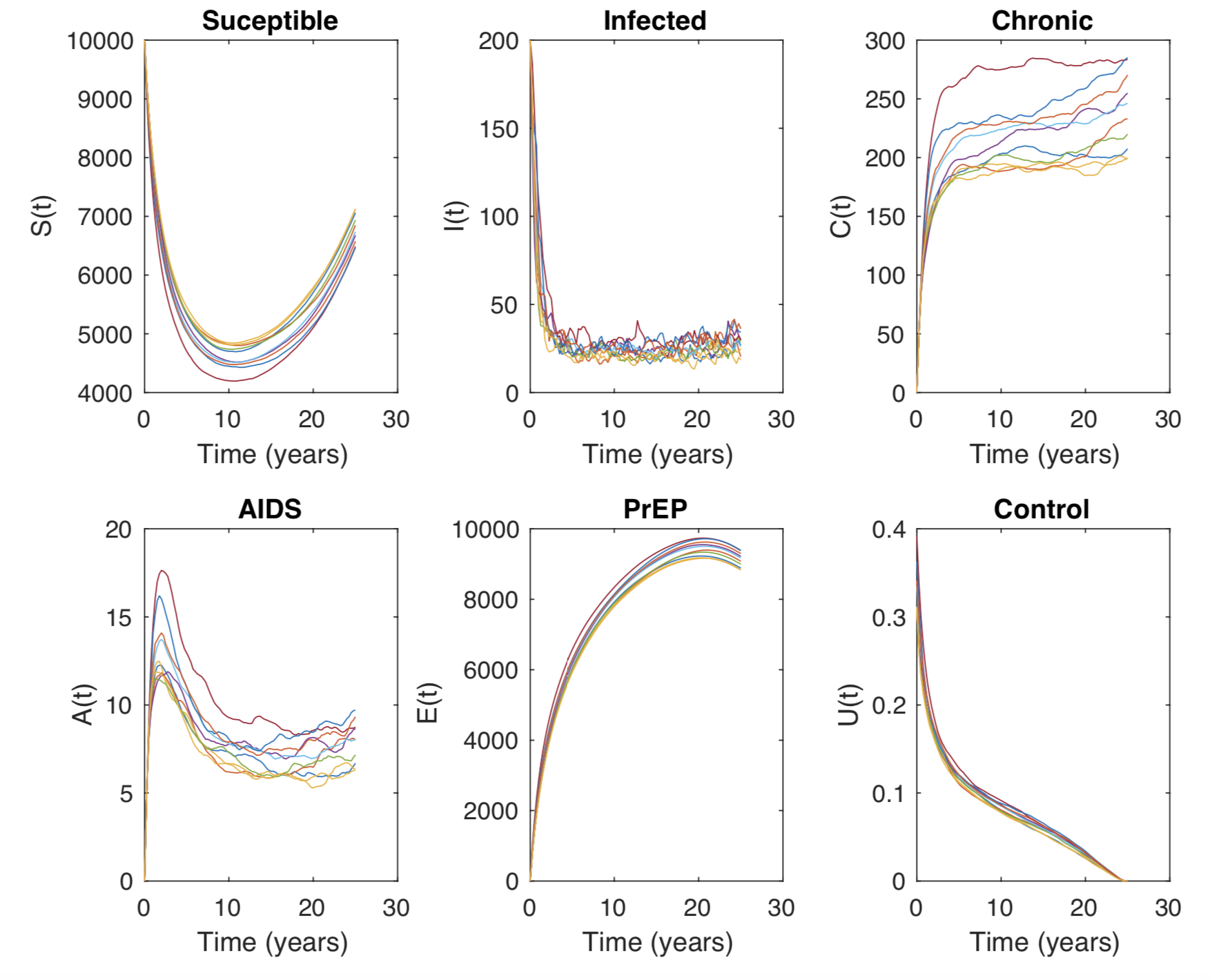}
\end{centering}
\caption{Plot of 10 paths of the solution of the stochastic optimal control PReP problem over 25 years with $\sigma =0.6/N$ and $N=30 000$.}
\label{fig: opt_control4}
\end{figure}

In Figure \ref{fig: opt_control4}, we have plotted the PReP model under the stochastic optimal control with $\sigma=0.6/N$. Hence, there is a larger weight on the noise terms in Figure \ref{fig: opt_control4} in comparison to that in Figure \ref{fig: opt_control3}.

By comparing Figure \ref{fig: opt_control} to Figure \ref{fig: opt_control3} and Figure \ref{fig: opt_control2} to Figure \ref{fig: opt_control4}, we see that the increase in $N$ leads to slightlydifferently shaped processes. By looking at the scale, we see that the optimal PReP treatment rate is significantly lower in the $N=30 000$ case than in the $N=10 200$ case. This is to be expected, since we have kept the initial values constant, and just increased the $N$, we have reduced the initial percentage of infected individuals in the population. Note also that the noise appears to affect the model more in the $N=30 000$ case of Figures \ref{fig: opt_control3} and \ref{fig: opt_control4} than in the $N=10 200$ case in Figures \ref{fig: opt_control} and \ref{fig: opt_control2}. 

\begin{figure}
\begin{centering}
\includegraphics[width=0.75\textwidth]{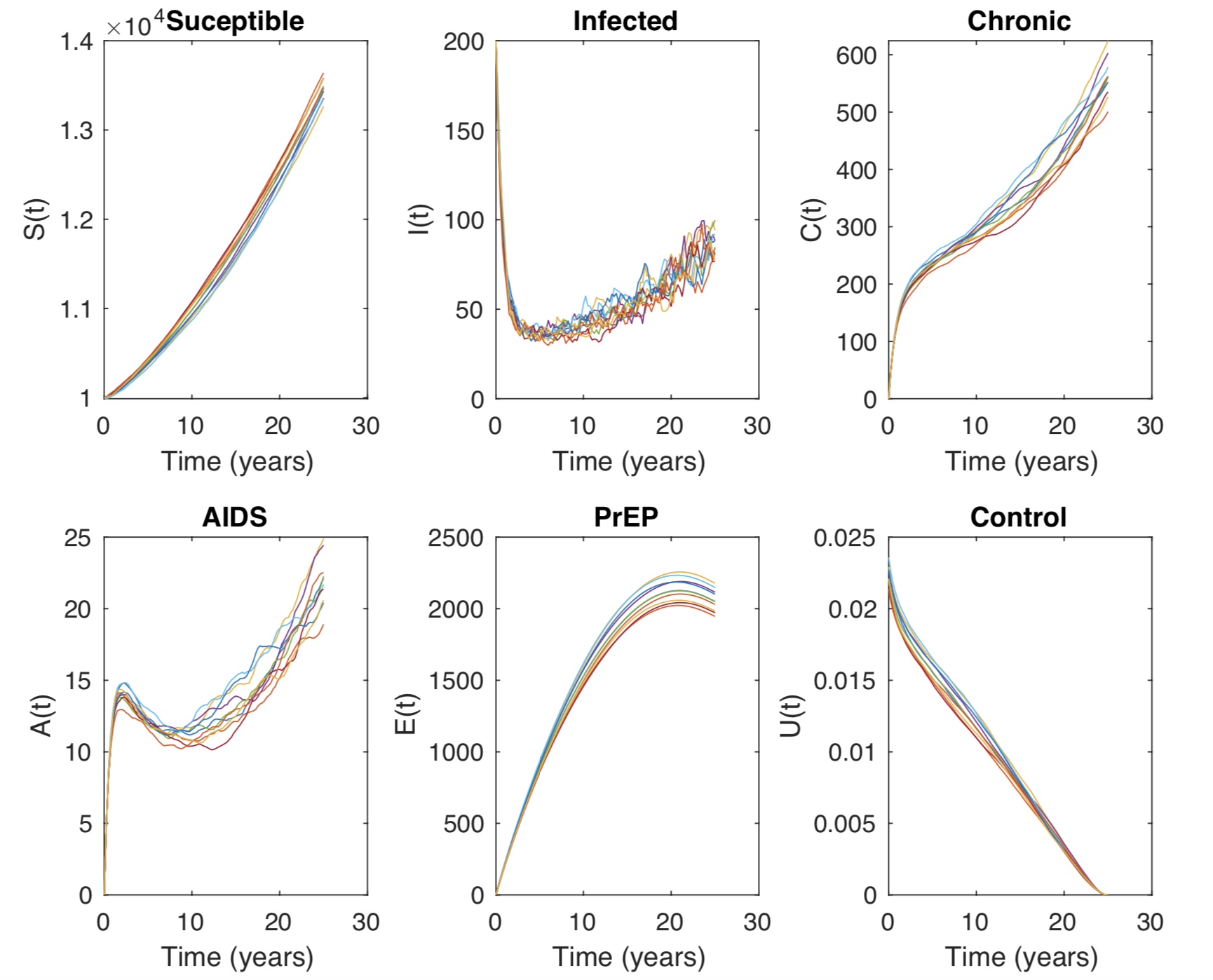}
\end{centering}
\caption{Plot of 10 paths of the solution of the stochastic optimal control PReP problem over 25 years with $\sigma = 0.2/N$ and $N=30 000$. Weights $w_1=0.2$, $w_2=0.3/N$.}
\label{fig: opt_control5}
\end{figure}

\begin{figure}
\begin{centering}
\includegraphics[width=0.75\textwidth]{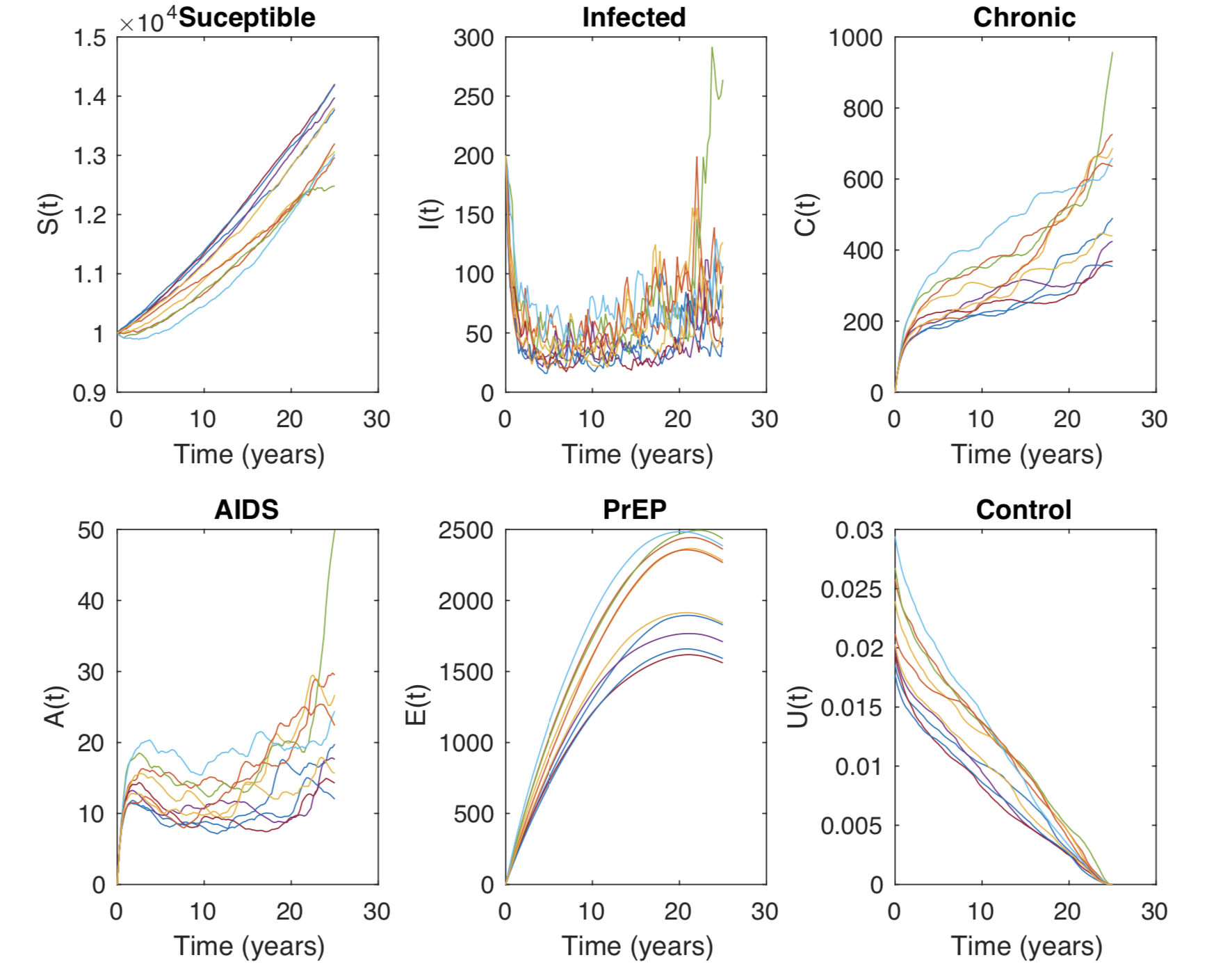}
\end{centering}
\caption{Plot of 10 paths of the solution of the stochastic optimal control PReP problem over 25 years with $\sigma =0.6/N$ and $N=30 000$. Weights $w_1=0.2$, $w_2=0.3/N$.}
\label{fig: opt_control6}
\end{figure}

In Figures \ref{fig: opt_control5} and \ref{fig: opt_control6}, we have kept everything from Figures \ref{fig: opt_control3} and \ref{fig: opt_control4} fixed, except the weights in the performance function. In Figures  \ref{fig: opt_control3} and \ref{fig: opt_control4}, we chose $w_1 = 20, w_2=0.3/N$. In Figures \ref{fig: opt_control5} and \ref{fig: opt_control6}, we instead choose $w_1 =0.2, w_2 =0.3/N.$

By looking at Figure \ref{fig: opt_control5} and Figure \ref{fig: opt_control6}, we see another effect of the type of performance function we have chosen: Since no weight is added to the number of infected individuals at the terminal time, it turns out to be optimal to let the number of infected individuals increase somewhat towards the end of the time period. This may seem counterintuitive, but is a consequence of our choice of performance function. To avoid this kind of optimal control, it would be better to consider a performance function of the form

\begin{equation}
\label{eq: performance_alt}
J(u) := E\left[\int_0^T [w_1 I(t) + w_2 u^2(t)]dt + w_3 I(T)\right],
\end{equation} 

where a weight is added to the number of infected individuals at the terminal time.

Another alternative is to consider the stochastic optimal control problem until infinite time:

\begin{equation}
\label{eq: performance_alt2}
J(u) := E\left[\int_0^{\infty} [w_1 I(t) + w_2 u^2(t)]dt\right].
\end{equation} 

These kinds of performance functions can be analysed by the same methods as in this paper, using the stochastic maximum principle and solving the corresponding SDE and adjoint BSDE numerically.

\section{Conclusions and future work}

The paper provides a study on stochastic optimal control for  the vaccine PReP, both in the unconstrained case and under two different kinds of constrains. The results show how stochastic optimal control theory can be applied for practical analysis in connection to control of infectious diseases.

Theoretical results were proven for the general case where the controlled PReP stochastic model for the spread of HIV may include jumps. An idea for future work is to prove existence and uniqueness results for the jump case of the stochastic model for the spread of HIV with PReP treatment. Also, conditions for the extinction and persistence of the disease would have to be established in the jump-case.

Further details about the numerical illustrations in Section  \ref{sec: numerical} will be provided in a forthcoming paper, as well as a more detailed mathematical analysis of the numerical methods. We would also like to numerically illustrate the stochastic optimal control results under the Type $I$ or Type $II$ constrains.

\bigskip

\textbf{Funding:} This research was partially supported by the  "Functional analysis and applications", Project 174007, MNTRS (Jasmina Djordjevi\'c) and "SCROLLER: A Stochastic ContROL approach to machine Learning with applications to Environmental Risk models", Project 299897 from the Norwegian Research Council (Kristina Rognlien Dahl).


\begin{thebibliography}{}

\bibitem{Baghery}
Baghery, F., \& {\O}ksendal, B., {\it A maximum principle for stochastic control with partial information}, Stochastic Analysis and Applications, 25(3),  (2007), 705-717.

\bibitem{Bismut}
Bismut, J. M., {\it Conjugate Convex Functions in Optimal Stochastic Control, Journal of Mathematical Analysis and Applications}, Vol. 44, (1973), 384-404.

\bibitem{Buckdahn}
Buckdahn, R., Djehiche, B. and Li, J., 
{\it A General Stochastic Maximum Principle for SDEs of Mean-field Type,Appl Math Optim }, (2011), 64:197-216 DOI 10.1007/s00245-011-9136-y

\bibitem{Campos}
C. Campos, C. J. Silva and D. F. M. Torres, \emph{Numerical optimal control of HIV transmission in Octave/MATLAB}, Mathematical and Computational Applications, 25, (2020), doi:10.3390/mca25010001.

\bibitem{DahlStokkereit}
K. R. Dahl and E. Stokkereit, {\it Stochastic maximum principle with Lagrange multipliers and optimal consumption with L{\'e}vy wage},  Afrika Matematika. Vol 27, (2016),   555-572.

\bibitem{DjordevicSilva}
J. Djordjevi\'c, Cristiana J. Silva,  {\it A stochastic analysis of the impact of fluctuations in the environment on pre-exposure prophylaxis for HIV infection}, Soft Computing, (2019), https://doi.org/10.1007/s00500-019-04611-1.


\bibitem{ElKarouiEtAl}
N.~El Karoui, S.~Hamadane, A.~Matousse, \newblock
{\it Backward Stochastic Differential Equations and Applications}, \newblock
\emph{Indifference Pricing: Theory and Applications} \newblock
Springer, Berlin Heidelberg, (2008),  267-320.

\bibitem{FOS}
N.~C.~Framstad, B.~{\O}ksendal, A.~Sulem, \newblock
\emph{Sufficient stochastic maximum principle for optimal control of jump diffusions and applications to finance}, \newblock
J. Opt. Theor. Appl., 121, (2007),  77-98.

\bibitem{GK}
Garcke, J.,  Kr\"oner, A., {\it Suboptimal feedback control of PDEs by solving HJB equations on adaptive sparse grids}, Journal of Scientific Computing, 70(1), (2017), 1-28.

\bibitem{KK}
Kalise, D., K.  Kunisch, {\it Polynomial Approximation of High-Dimensional Hamilton--Jacobi--Bellman Equations and Applications to Feedback Control of Semilinear Parabolic PDEs},  SIAM Journal on Scientific Computing, 40(2),  (2018), A629-A652.

\bibitem{Kushner}
Kushner, H. J., {\it Necessary Conditions for Continuous Parameter Stochastic Optimization Problems}, SIAM Journal on Control, Vol. 10, (1972),  550-565.

\bibitem{Oksendal_delay}
{\O}ksendal, B and Sulem, A.. In J.M. Menaldi, E. Rofman and A. Sulem (editors),
{\it A maximum principle for optimal control of stochastic systems with delay, with applications to finance}, 
: Optimal Control and Partial Differential Equations - Innovations and Applications, IOS Press, Amsterdam 2000.

\bibitem{Nichols}
Nichols, B.E., Boucher, C.A., van der Valk M., Rijnders B. J., van de Vijver D.A., \emph{Cost-effectiveness analysis of pre-exposure prophylaxis for 691 HIV-1 prevention in the Netherlands: a mathematical modelling 692 study}, Lancet Infect Dis (2016)  16(12):1423-1429.

\bibitem{Oksendal}
{\O}ksendal, B., {\it Stochastic Differential Equation}, Springer, Berlin Heidelberg, 6th ed., 2007.

 \bibitem{OS-Levy}
{\O}ksendal, B. and  Sulem, A., \newblock
\emph{Applied Stochastic Control of Jump Diffusions}, \newblock
2. ed., Springer, Berlin Heidelberg, (2007).

\bibitem{OksendalSulem}
{\O}ksendal, B.  and Sulem,  A., {\it Risk minimization in financial markets modeled by It{\^o}-L{\'e}vy processes}, Afrika Matematika, (2015),  vol 26,  939-979

\bibitem{Perelson}
Perelson, A. S. et al,  {\it Decay characteristics of HIV-1-infected compartments during combination therapy}, Nature (1997) 387:188-191.

\bibitem{Peng}
Peng, S.,  {\it A General Stochastic Maximum Principle for Optimal Control Problems}, SIAM J. Control Optim., vol. 28, (1990), 966-979.

\bibitem{SilvaTorres2017}
Silva, C. J. and Torres, D. F. M., {\it  A SICA compartmental model in epidemiology with application to HIV/AIDS in Cape Verde.},  Ecol Complex 708,  (2017), 30:70-75.


\bibitem{SilvaTorres}
Silva, C. J. and Torres, D. F. M., {\it Modeling and optimal control of HIV/AIDS prevention through PReP}, Discrete and Continuous Dynamical Systems Series S, vol. 11, no. 1, (2018), 119-141.

\bibitem{Sharomi}
Sharomi, O., Podder, C. and Gumel, AB., {\it Mathematical analysis of the  transmission dynamics of HIV/TB co-infection in the presence of treatment}, Math Biosci Eng (2008) 5:145-174.


\bibitem{TangLi}
Tang, S. and Li, X., \newblock
\emph{Necessary conditions for optimal control of stochastic systems with random jumps}, \newblock
SIAM Journal of Control and Optimization, 5, (1994),1447-1475.


\bibitem{Zwahlen}
Zwahlen, M. and Egger, M., {\it Progression and mortality of untreated 732 HIV-positive individuals living in resource-limited settings}, update  of literature review and evidence synthesis. Report on UNAIDS 734 obligation no HQ/05/422204, (2006).

\end{thebibliography}
\end{document}